\documentclass{puthesis}

\author{Vyacheslav Rychkov}
\title{Estimates for Oscillatory Integral Operators}
\departmentprefix{Department of}
\department{Mathematics}

\abstract{
This thesis is devoted to asymptotic norm estimates for oscillatory integral
operators acting on the $L^2$ space of functions of one real variable.
The operators in question have compact support and an oscillatory
kernel of the form $\exp(i\gl S(x,y))$, where $S(x,y)$ is a smooth real
phase function, and $\gl$ is a large real number. 

I study how the norm of the operator decays as $\gl$ goes to
infinity, and how the rate of this decay can be determined from
the properties of the phase function $S(x,y)$.

For $\cin$ phase functions I prove results formulated in terms of the
Newton polygon of $S(x,y)$, improving previously known estimates. 
My estimates are best possible or differ from the best possible ones by
at most a power of $\log\gl$.

I use two different methods. The first method is based on the
geometric analysis of the zero set of the Hessian $\sxy$ using the
Puiseux decompositions. The second method is based on a stopping time argument. 

}
\acknowledgements{

I would like to thank most sincerely my mathematics advisor Elias Stein 
for believing in me, for teaching me, for giving me support,
understanding, and inspiration during all my studies at Princeton.
If it were not for Eli, I probably would have not come to Princeton,
and my life would have taken quite a different path. As for this
thesis, it would have never been written.

I would like to express my deep gratitude to my physics advisor 
Alexander Polyakov for changing my perception of science forever, and  
for the ideas, knowledge, and wisdom he shared with me so generously and 
patiently. Without A.M. my conversion to physics would 
never have happened.

I am infinitely indebted to my former advisors Oleg V. Besov of Moscow, Gennady
A. Kalyabin of Samara, and Hans Triebel of Jena, who taught me how to do mathematics 
in Russia and Germany before I came to Princeton. 

I talked about the results of my thesis in the Analysis Seminar at
the University of Wisconsin at Madison in the Fall of 2001. I would
like to thank
Andreas Seeger, Stephen Wainger, and other analysts at UW for their
kind invitation, hospitality, and interesting discussions.

I would like to thank my colleague and friend Igor Rodnianski for reading
this thesis in a highly professional and timely fashion, and for the
remarks he made, questions he asked, and misprints he caught.

I would like to thank our graduate coordinator Jill LeClair, whose
care and expertise magically transformed the everyday 
graduate student life at the Mathematics Department from chaos to order. 
 
Finally, I would like to thank my dear friends Shura Schekochihin, Vamsi
Madhav, Jorge Silva, Jeroen van Dongen, Felix Vardy, Harald Helfgott, 
Etienne De Poortere, Maria Eug\^enia Dias Cruz, Andrei Efanov,
Vadim Kaloshin, Olesya Grischenko, Sasha Baitin, Sasha
 Ulianova, Nadia Zakamskaya, Lenia Koralov, Dima Sarkisov, Edith
Elkind, Toufic Suidan, Seth Patinkin,
Gergely Harcos, Martin Niepel, Eaman Eftekhary, Hadi Jorati, Yannis
Papadoyannakis, Mihalis Dafermos
for all the 
invaluable help, learned advice, worthy example, immeasurable joy,
enlightened simplicity, and perverted sophistication 
they gave me. I cannot imagine my life in Princeton
without them.

Above all, I thank my parents Sergei and Tatiana,
to whom I owe everything.
\begin{flushright}
Princeton, April 4, 2002
\end{flushright}

}
\dedication{To the memory of my grandfathers,\\Vyacheslav Vasilievich
Chelpanov\\and Pavel Pavlovich Rychkov}

\usepackage{amssymb}
\usepackage{amsthm}
\usepackage[dvips]{graphicx}

\newcommand{\NC}{\newcommand}

\newcommand{\gf}{\varphi} \newcommand{\eps}{\varepsilon}
\newcommand{\gl}{\lambda} \newcommand{\ga}{\alpha}
\newcommand{\gb}{\beta}  \newcommand{\gd}{\delta}
\newcommand{\gp}{\psi}   
 
\renewcommand{\gg}{\gamma}  \newcommand{\gs}{\sigma}

 \newcommand{\GP}{\Psi}
\newcommand{\GF}{\Phi}   
\newcommand{\GG}{\Gamma} \newcommand{\GD}{\Delta}
\newcommand{\GS}{\Sigma}

\newcommand{\<}{\le} 
\renewcommand{\>}{\ge}
\NC{\gsim}{\gtrsim}
\NC{\lsim}{\lesssim}

\renewcommand{\cal}{\mathcal}
\renewcommand{\bar}{\overline}
\newcommand{\dbarchi}{\overline{\overline\chi}}
\newcommand{\supp}{\mathop{{\rm supp}}}
\renewcommand{\tilde}{\widetilde}

\newcommand{\bb}{ \mathbb}
\newcommand{\rr}{{\bb R}}
\newcommand{\cc}{{\bb C}}
\newcommand{\nat}{{\bb N}}

\newcommand{\z}{{\bb Z}}
\newcommand{\cin}{C^\infty}

\newcommand{\sumi}{\sum^\infty}
\NC{\rjk}{R_{jk}}
\NC{\tjk}{T_{jk}}
\NC{\Picture}[1]{
\begin{center}
\scalebox{0.8}{\includegraphics{#1}}
\end{center}}

\newcommand{\rf}[1]{{\rm(\ref{#1})}}
\newcommand{\dis}{\displaystyle}

\newcommand{\sxy}{S''_{xy}}
\newcommand{\ltwor}{L^2(\rr)}
\newcommand{\nt}{\|T\|}
\NC{\tl}{T_\lambda}
\NC{\ntl}{\|\tl\|}
\newcommand{\zp}{\z_+}
\newcommand{\rrp}{\rr_+}

\renewcommand{\)}{)\!)}
\newcommand{\cinxy}{\cin\(x,y\)}
\newcommand{\tup}[1]{\textup{(#1)}}
\newcommand{\D}{\partial}
\renewcommand{\Re}{{\rm Re\,}}
\renewcommand{\Im}{{\rm Im\,}}
\newcommand{\dist}{{\rm dist}}

\renewcommand{\qedsymbol}{$\blacksquare$}

\NC{\beq}{\begin{equation}}
\NC{\beql}[1]{\beq\label{#1}}
\NC{\eeq}{\end{equation}}

\NC{\begi}{\begin{itemize}}
\NC{\ei}{\end{itemize}}

\newtheorem{Th}{Theorem}[chapter] 
\newtheorem{Lem}[Th]{Lemma} 

\newtheorem{Prop}[Th]{Proposition} 
\theoremstyle{definition}

\theoremstyle{remark}

\begin{document}

\chapter{Introduction}

\section{Formulation of the problem}

My thesis studies asymptotic norm estimates for oscillatory integral
operators acting on the $L^2$ space of functions of one real variable.

More precisely, I fix a real $\cin$ function $S(x,y)$ (called
{\it phase function}) and consider a one-parameter family of operators of the form
\beq
\label{operator}
T_\lambda f(x)=\int_{-\infty}^{\infty}e^{i\gl S(x,y)}\chi(x,y)f(y)\,dy\qquad(\lambda\in\rr). 
\eeq
Here $\chi(x,y)$ is an unimportant $\cin$ cut-off function
compactly supported in a small neighborhood of the origin in $\rr^2$. 

The operators $T_\lambda$ act on $\ltwor$, and it is generally to be expected
that for $\gl\to\infty$ the norm $\|T_\lambda\|$ will
decay. Typically, we will have 
\beq
\label{decay}
\|T_\lambda\|\le C \lambda^{-\delta}\qquad(\lambda\to\infty)
\eeq
My thesis studies how the decay rate $\delta$ depends on
the properties of the phase function $S(x,y)$.

\section{History and motivation}

H\"ormander \cite{Hor} proved that if $\sxy\ne0$ on the support of
$\chi$, then \rf{decay} is true with $\gd=1/2$, and this is best
possible.

Typical example falling under the scope of this result is $S(x,y)=xy$.
For this choice of the phase functions $T_\lambda$ is a rescaled and
cut-off version of the Fourier transform.

Some problems of harmonic analysis naturally lead to more general
phase functions which do not necessarily satisfy H\"ormander's
condition $\sxy\ne0$. 

For instance, such phase functions arise in studying smoothing properties of generalized
Radon transform associated with families of curves having various
geometric degeneracies. The most direct connection exists between 
decay norm estimates for $T_\lambda$ and smoothing properties of the generalized Radon
transform in the plane defined by 
$$
{\cal R}g(t,x)=\int_{-\infty}^\infty g(t+S(x,y),y)\chi(x,y)\,dy.
$$
Namely the decay estimate (\ref{decay}) implies that ${\cal R}$ is a smoothing operator of order $\delta$,
that is, it acts from the Sobolev space $H^s(\rr^2)$ into
$H^{s+\delta}(\rr^2)$ for any $s$ (see Phong and Stein \cite{PSstop}).

\section{Known results for degenerate phase function}
 
As I mentioned above, it is desirable to be able to determine
the optimal exponent $\delta$ in (\ref{decay}) for phase
functions with vanishing $\sxy$, which are called {\it degenerate}.

My thesis addresses this problem in its local aspect.
That is, I concentrate on the properties of $S(x,y)$ near the origin, 
and use the freedom to choose the support 
of the cut-off function $\chi(x,y)$ as small as I want.
The typical form of results that I will state is going to be ``There exists a small
neighborhood of the origin $U$ such that if $\supp\chi\subset U$ 
, then \ldots''.
 
As it was realized by Phong and Stein \cite{PS}, the optimal exponent
$\delta$ depends on the local properties of the phase function
$S(x,y)$ at the origin via the Newton polygon of $\sxy$. 

The Newton polygon of $\sxy$ is defined as follows. In the positive
quadrant of the plane mark all the points with integer coordinates
$(p,q)$ such that the partial derivative
$$
\left(\frac\partial{\partial x}\right)^p\left(\frac\partial{\partial y}\right)^q \sxy(0,0)\neq 0.
$$
After that, take the marked point for which $p$ is minimal and add the
vertical ray emanating from this point upward. Also, take the marked
point for which $q$ is minimal and add the horizontal ray emanating
from this point leftward. The Newton polygon is the convex hull of the
set consisting of all marked points and two added rays (see Fig. 1).

\Picture{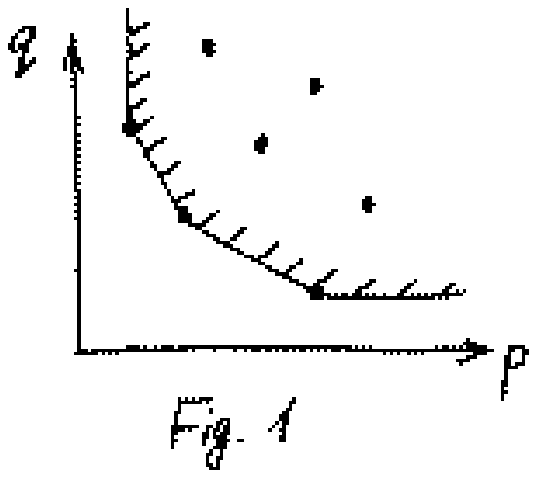}

Assume that the Newton polygon of $\sxy$ is not empty, which means
that not all partial derivatives of $\sxy$ vanish at the origin.
Denote by $t_0$ the parameter of intersection of the line $n_1=n_2=t$
with the boundary of the Newton polygon. The number
$$
\GD=\frac 1{2t_0+2}
$$
is called the {\it 
Newton decay rate} of $S(x,y)$ (this definition differs by a factor of
1/2 from \cite{PS} and \cite{R}). 

Known norm estimates for $T_\lambda$ relevant for my
work are the following:
\begin{itemize}
\item (Phong and Stein \cite{PS}) 
Lower estimate $\ntl\>C\gl^{-\GD}$. 
\item (implicitly in Seeger \cite{See1},\cite{See2}) Almost sharp upper estimate
$\ntl\<C_\eps \gl^{-\GD+\eps}$ for any $\eps>0$.
\item (Phong and Stein \cite{PS}) Sharp upper estimate $\ntl\<C
\gl^{-\GD}$ under the additional assumption that $S(x,y)$ is real analytic.
\end{itemize}
The upper estimates are true provided the support of $\chi(x,y)$ is small
enough. The lower estimate is true for $\chi(0,0)\neq 0$. In all three
estimates $\gl\to\infty$.

\section{Main result of the thesis}

The purpose of my thesis is to show that the sharp estimate proven by
Phong and Stein in the real analytic case continues to hold in the
$\cin$ case without loss of $\eps$. There will be one possible exception, when one
loses at most a power of log. 

Consider the formal Taylor series of $\sxy$ at the origin
$$
\sxy(x,y)\sim\sum_{p,q} c_{pq} x^p y^q,\qquad
c_{pq}=\frac1{p!q!}\D^p_x\D^q_y \sxy(0,0).
$$
I say that $\sxy$ is \emph{exceptionally degenerate}, if this series
can be factored in the ring of formal power series 
$\rr[[x,y]]$ into the product
\[
U(x,y)(y-f(x))^N,
\]
where 
\begin{itemize}
\item
$N\>2$,  
\item the series $f(x)\in\rr[[x]]$ is of the form $f(x)=cx+\ldots$
with $c\ne0$, and
\item
the series $U(x,y)\in\rr[[x,y]]$ is invertible, that is its zeroth
order term is nonzero. 
\end{itemize}
Note that $\GD=\frac1{N+2}$ for such a phase function. 

The main result of the thesis is the following
\begin{Th} 
\label{main}
There exists a small
neighborhood
of the origin $V$ such that\\
\textup{(a)}
If $\sxy$ is not exceptionally degenerate, and $\supp\chi\subset V$, 
then 
$$\ntl\<C\gl^{-\GD}\qquad(\gl\to\infty).$$
\textup{(b)} If $\sxy$ is exceptionally degenerate, and $\supp\chi\subset V$,
then 
\begin{equation}
\label{main1}
\ntl\<C\gl^{-\frac1{N+2}}(\log\gl)^{\frac{2N-1}{N+2}}\qquad(\gl\to\infty).
\end{equation}
\end{Th}

This theorem was proved in my paper \cite{R}, modulo an inessential
improvement of the power of $\log\gl$ in \rf{main1}. In Chapter 5 I will
show that in fact for $N=2$ estimate \rf{main1} can be improved to the
sharp one $\ntl\<C\gl^{-1/4}$. I do not know if a similar improvement
is possible for $N\>3$.

\chapter{Real analytic case}

In this chapter I give an argument for the upper estimate of
Phong and Stein in the real analytic case. This argument is somewhat simpler
than the original proof. The main purpose is to familiarize the reader with
the ideas and the technology, which will later be partially recycled in
the proof of Theorem 1.1. I will explain what the main difficulty is
going to be in generalizing to the $\cin$ case. To begin with, I
review the lower bound. In this chapter no new results are proved.

\section {Lower bound}

The proof of the lower bound $\ntl\gtrsim \gl^{-\GD}$ is obtained by
looking at the regions of the $(x,y)$ plane where $S(x,y)\approx const
\sim 1/\gl$ and restricting the operator to those regions.

\subsection{Typical example}

I consider the example of the polynomial phase
function whose
Newton polygon has only two corner points:
\beq
\label{ex}
S(x,y)=c_1 x^{A_1}y^{B_1}+c_2  x^{A_2}y^{B_2}\qquad(B_2>B_1).
\eeq
I also assume in this example that the segment joining points
$(A_i,B_i)$
intersects the bisectrix of the $(A,B)$ plane (Fig. 2).
This example captures the main idea of the proof in the general case.

\Picture{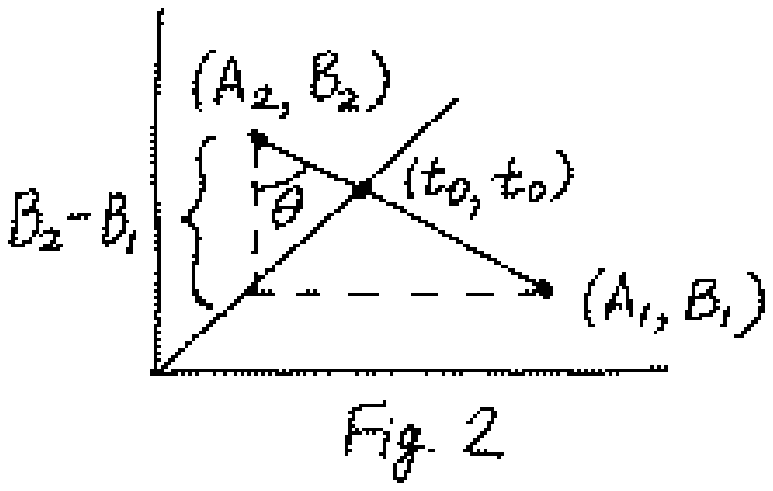}

We can expand $S(x,y)$ as follows:
$$
S(x,y)=Cx^{A_2}y^{B_1}\prod_{i=1}^{B_2-B_1}(y-\theta_ix^\gamma),\qquad
\gamma=\frac{A_1-A_2}{B_2-B_1},\quad\theta_i\in\cc.
$$
Notice that the exponent $\gamma=\tan \theta$ in Fig. 2, while the number of
branches $B_2-B_1$ is equal to the height of the triangle. The reader
will see later that these features naturally extend to general phase
functions.

Now notice that for $y=C x^\gg$ for generic $C$ different from all
$\theta_i$ I have
\beq
\label{xg}
S(x,y)=const. x^{A_2+\gg B_1+\gg(B_2-B_1)}\sim x^{A_2+\gg B_2}.
\eeq 

 From the equation of the
straight line passing through $(A_i,B_i)$,
$$
\frac{x-A_2}{y-B_2}=\frac{A_1 -A_2}{B_1-B_2}=-\gg,
$$
putting $x=y=t_0$ I find
$$
t_0=\frac{A_2+\gg B_2}{1+\gg}.
$$

I am looking for a region of the $(x,y)$-plane where $y\sim x^\gg$ and $S(x,y)\sim
1/\gl$. According to \rf{xg}, this happens for 
$$
x\sim\gl^{-\frac 1{t_0(1+\gg)}},\qquad y\sim\gl^{-\frac
\gg{t_0(1+\gg)}}.
$$
It follows that I can find a rectangle $R$ of size $\gd_x\times\gd_y$ with
sides parellel to the coordinate axes with
$$
\gd_x\sim\gl^{-\frac 1{t_0(1+\gg)}},\qquad \gd_y\sim\gl^{-\frac
\gg{t_0(1+\gg)}},
$$
and such that $\Re e^{i\gl S(x,y)}>c>0$ on $R$.

It remains to invoke the following straighforward

\begin{Lem} Suppose the kernel $K(x,y)$ of an integral operator
\beq
\label{opT}
Tf(x)=\int K(x,y) f(y)dy
\eeq
is such that $\Re K(x,y)>c>0$ on a rectangle of size
$\gd_x\times\gd_y$. Then the norm of the operator on $\ltwor$
satisfies $\nt\> c(\gd_x\gd_y)^{1/2}$.
\end{Lem}

Using the lemma, I get 
\beq
\label{low}
\ntl\gtrsim (\gd_x\gd_y)^{1/2}\sim\gl^{-1/(2t_0)}.
\eeq
This is the correct answer, since I was looking at the Newton polygon
of $S$, which differs from the Newton polygon of $\sxy$ by a shift by
vector $(1,1)$.

\subsection{General case}

The general case turns out to be very similar to the example I have
just considered. The argument does not use real analyticity and 
works generally in the $\cin$ case.

I have an asymptotic expansion
\beq
\label{asexp}
S(x,y)\sim\sum C_{AB}x^Ay^B,
\eeq
where $(A,B)$ runs through points inside the Newton polygon of
$S(x,y)$. I consider the edge of the Newton polygon intersecting the
bisectrix (the {\it main edge}) and look at the region $y\sim
x^\gamma$, where $\gamma$ is the exponent associated with this
edge.

A simple check shows that in this region the terms in \rf{asexp} coming from all the
edges but the main one and from the inside of the Newton polygon are
subleading.
As a result, \rf{xg} is true as before, and \rf{low} follows. \qedsymbol

\section{Upper bound}

As I have just shown, the lower bound follows by restricting the
operator to a rectangle where the phase function is effectively
constant,
so that the oscillatory behavior is suppressed. The upper bound, to
which we proceed, is much trickier. It requires a decomposition of the
$(x,y)$ plane into much bigger rectangles, on which the phase function
does oscillate, in a controlled fashion.

\subsection{Elementary tools}

The following 3 elementary results are
needed in the proof of the upper bound. In a sense, in most cases
you just need to find the right combination of the tools which works,
and do the algebra corectly.

All rectangles below are assumed to have sides parallel to the
coordinate axes.

\begin{Lem}[Size estimate]
\label{size}
Let $T$ be an integral operator of the form \rf{opT}. Assume that the
kernel $K(x,y)$ is supported in a rectangle of size $\gd_x\times\gd_y$
and bounded: $|K|\<1$. Then $T$ is bounded on $\ltwor$ with the norm
$\nt\<(\gd_x \gd_y)^{1/2}$. 
\end{Lem}

\begin{Lem}[Oscillatory estimate] 
\label{oscest} 
Let $\tl$ be an oscillatory integral operator of the form
\rf{operator}. Assume that\\
\tup{1} $\chi(x,y)$ is supported in a
rectangle $R$ of size $\gd_x\times\gd_y$,\\
\tup{2} $|\D_y^n\chi|\<C/\gd_y^{n}$  in  $R$ for $n=0,1,2$,\\
\tup{3} $|\sxy|\>\mu>0$  in  $R$,\\
\tup{4} $|\D_y^n\sxy|\<C\mu/\gd_y^n$  in  $R$ for $n=0,1,2$.\\
Then $\ntl\<const(\gl\mu)^{-1/2}$ with $const$ depending only on $C$.
\end{Lem}

\begin{Lem}[Almost orthogonality]
\label{aa}
Let $\{R_j\}$ be a family of rectangles and $\{T_j\}$ be a family of integral
operators of the form \rf{opT} such that\\ 
\tup{1} the kernel of $T_j$ is supported in $R_j$,\\
\tup{2} the family $\{R_j\}$ is almost orthogonal in the sense that for
every rectangle $R_j$ the number of rectangles whose horisontal
or vertical pojections intersect those of $R_j$ is bounded by a
constant $C$.\\
\tup{3} $T_j$ are bounded on $\ltwor$ with $\|T_j\|\<A$ independent of
$j$.\\
Then $T=\sum_j T_j$ is bounded on $\ltwor$ with norm $\nt\<const.A$, where
$const$ depends on $C$ only.  
\end{Lem}

The proof of the size estimate is straightforward (consider
$(Tf,g)$). The oscillatory estimate 
is a variant of the \emph{Operator van der Corput lemma} of Phong
and
Stein \cite{PS}. The lemma is proved by a standard $TT^*$
argument. The assumptions made are enough to show, integrating
by parts twice, that the kernel of $TT^*$ has the bound
\[
K(x_1,x_2)\<C\frac{\gd_y}{1+\gl^2\mu^2\gd_y^2|x_1-x_2|^2},
\]
which implies the necessary norm estimate. We omit the details.
The almost orthogonality is a trivial consequence of the Cotlar-Stein lemma.

\subsection{Example}

Once again, I start with an example. This time I take the $\sxy$
rather then $S$ in the form of \rf{ex}
$$
 \sxy(x,y)=c_1 x^{A_1}y^{B_1}+c_2  x^{A_2}y^{B_2}\qquad(B_2>B_1).
$$
While I looked at the phase $S(x,y)$ when proving the lower bound,
it is its second mixed derivative $\sxy$ which is important for the
upper bound.

I also assume again that the segment joining points
$(A_i,B_i)$ is the main edge of the Newton polygon of $\sxy$, that is it
intersects the bisectrix of the $(A,B)$ plane (Fig. 2).

The proof starts by taking the dyadic decomposition of the $(x,y)$
plane into rectangles $R_{jk}$ of size $2^{-j}\times 2^{-k}$.  
I take a suitable smooth partition of unity fitted to this family of
rectangles,
and use it to localize the operator $\tl$ to $R_{jk}$, that is, I
consider the operators
$$
T_{jk} f(x)=\int e^{i\gl
S(x,y)}\chi_{jk}(x,y)\chi(x,y)f(y)dy,\qquad\sum T_{jk}=\tl,
$$
where $\supp\chi_{jk}$ is contained in the doubled rectangle
$R_{jk}^*$ (Fig. 3).

\Picture{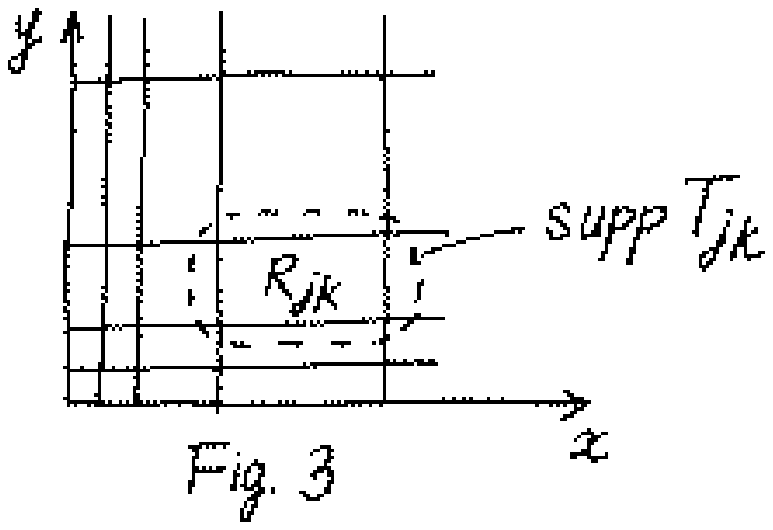}

I again look at the expansion
\beq\label{exp1}
\sxy(x,y)=Cx^{A_2}y^{B_1}\prod_{i=1}^{B_2-B_1}(y-\theta_ix^\gamma).
\eeq
I notice that there are 3 important regions of parameters
$(j,k)$. Case I: $k\<\gamma j-K$, where $K$ is a large constant. 
Case II: $k\>\gg j+K$. Case III: $|k-\gg j|\<K$. These regions corespond to
rectangles $R_{jk}$ lying respectively well above, well below, and
around the curve $y=x^\gg$ (Fig. 4).

\Picture{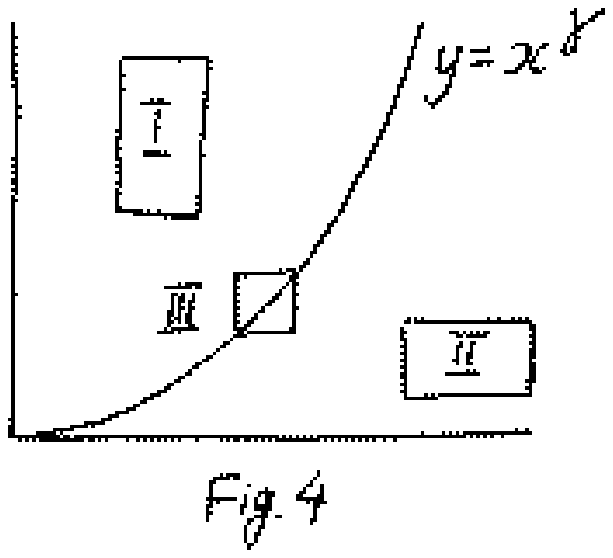}

 {\bf Case I.} By the size estimate (Lemma \ref{size}) I know the individual bounds
\beq\label{size1}
\|T_{jk}\|\lesssim 2^{-(j+k)/2}.
\eeq
To apply the oscillatory estimate, I need to know how large $\sxy$ is
on $R_{jk}$. It is easy to see from \rf{exp1} that provided $K$ is
chosen large enough, I have the following estimates on $R^*_{jk}$
\beq\label{osc11}
|\sxy|\sim 2^{-jA_2-kB_2}.
\eeq
In this particular situation it is easy to check that 
the remaining conditions of Lemma \ref{oscest} are satisfied, so that I
conclude
\beq\label{osc1}
\|T_{jk}\|\lesssim \gl^{-1/2}  2^{(jA_2+kB_2)/2}\qquad(k\<\gg j-K) .
\eeq

The most natural thing to do is to take a geometric mean
$a^{1-\theta}b^\theta$ of estimates
\rf{size1} and \rf{osc1}, choosing $\theta=2\GD$, 
so that the resulting estimate will have the desired $\gl$
behavior $\gl^{-\GD}$. Remember that 
$$
\GD=\frac 1{2t_0+2}=\frac{1+\gg}{2(A_2+1)+2\gg(B_2+1)}.
$$
Doing the algebra, I obtain the following estimate
\beq\label{ind}
\|T_{jk}\|\lesssim
\gl^{-\GD}2^{D(k-\gg j)},\qquad D=\frac{B_2-A_2}{2(A_2+1)+2\gg(B_2+1)}.
\eeq

Now it's time to invoke almost orthogonality.
I split all the Case I rectangles into families indexed by
a natural number $r$, putting into the $r$-th family all $R_{jk}$ such
that
$$
k-[\gg j]=-r.
$$
For each $r$, such a family is almost orthogonal, and so it follows
from \rf{ind} by Lemma \ref{aa} that 
\beq\label{estr}
\left\|\sum_{k-[\gg j]=-r}T_{jk}\right\|\lesssim\gl^{-\GD}2^{-Dr}\qquad(r\>K).    
\eeq

Notice that in general I have $B_2\>A_2$ (Fig. 2). Assume for the
moment that $B_2>A_2$. In this case $D>0$, and I can sum estimate
\rf{estr} over $r$ from $K$ to infinity, to get
\beq\label{estrr}
\left\|\sum_{k-\gg j\<-K}T_{jk}\right\|\lesssim\gl^{-\GD},
\eeq
which is the required estimate.

If $B_2=A_2$, I am going to start again from estimates \rf{size1} and
\rf{osc1}
and use a completely different splitting into almost orthogonal
families. In this case, I will put into the $r$-th family all $R_{jk}$
such that 
$$
j+k=r.
$$
By almost orthogonality, it follows from \rf{size1} and
\rf{osc1} that
\beql{jkr}
\left\|\sum_{j+k=r,\;k-\gg j\<-K}T_{jk}\right\|\lesssim\min(2^{-r/2},\gl^{-1/2}2^{rA_2/2}).
\eeq
The decreasing and increasing progressions under the minimum sign
balance for
$$
r=r_*=\frac{\log_2\gl}{A_2+1}\pm const.
$$
Summing \rf{jkr} in $r$, I get
\beq\label{a=b}
\left\|\sum_{k-\gg j\<-K}T_{jk}\right\|\lsim 2^{-r_*/2}\lesssim\gl^{-1/(2A_2+2)},
\eeq
which is exactly what is required.

 {\bf Case II.} It comes as no suprise that this case is going to be
absolutely similar to Case I. The size estimate \rf{size1} stays the
same, and the appropriate oscillatory estimate obtained analogously to
\rf{osc1} comes out to be
$$
\|T_{jk}\|\lesssim \gl^{-1/2}  2^{(jA_1+kB_1)/2}\qquad(k\>\gg j+K).
$$

If $B_1=A_1$, I split into almost orthogonal families $j+k=r$ and arrive
at the analogue of \rf{a=b}. 

If $B_1<A_1$ (notice that always $B_1\<A_1$), I do the same manipulation 
which led to \rf{ind}, and get
$$
\|T_{jk}\|\lesssim
\gl^{-\GD}2^{D'(k-\gg j)},\qquad D'=\frac{B_1-A_1}{2(A_1+1)+2\gg(B_1+1)}.
$$
Notice that now $D'<0$, which is exactly compatible with having to sum
over $k-\gg j\>K$. I split into the almost orthogonal fiamilies $k-\gg
j=r\>K$, and get the analogues of \rf{estr} and \rf{estrr}. Case
closed.

{\bf Case III.} Here I cannot get any reliable estimates on $\sxy$ on the
whole rectangle $R_{jk}$. Because of this, a further decomposition is
required. However, the present example is a bit too special to demonstrate
the method. I will deal with this situation in a more general setting 
in the next section. 

\subsection{General case}

Now I am going to consider the general case of real analytic $\sxy$. 
The basis of my consideration is going to be the following far-reaching generalization
of \rf{exp1} known as {\it Puiseux theorem}. This result is basically
well known (see \cite{PS}). 

\Picture{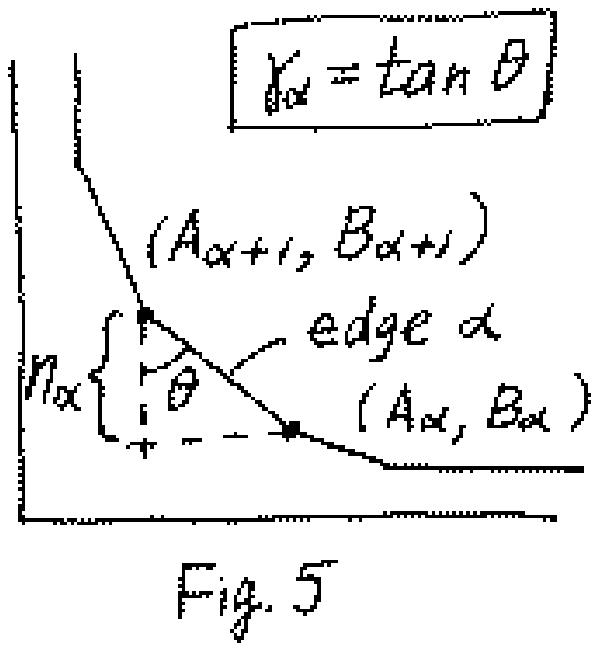}

First I look at the
Newton polygon of $\sxy$. In general, the polygon is going to have
some number of finite edges and two infinite edges (Fig. 5). With each finite
edge $\alpha$ joining points $(A_\ga,B_\ga)$ and
$(A_{\ga+1},B_{\ga+1})$, $B_{\ga+1}>B_\ga$, I associate numbers $\gg_\ga>0$ and
$n_\ga\in\nat$, where
$$
\gg_\ga=\frac{A_\ga-A_{\ga+1}}{B_{\ga+1}-B_\ga},\qquad n_\alpha=B_{\ga+1}-B_\ga.
$$
Let also $A$ and $B$ be the $x$ and $y$ coordinates of the infinite edges.
Then the claim of the Puiseux theorem is that in a neighborhood of the origin there exists a factorization
\beq\label{Pu}
\sxy=U(x,y)x^Ay^B\prod_\ga\prod_{i=1}^{n_\ga}(y-Y_{\ga i}(x)),
\eeq
where $Y_{\ga_i}(x)$ are convergent fractional power series with
fractionality at most $x^{1/n!}$, $n=B+\sum n_\ga$, whose
expansion starts with
$$
Y_{\ga i}(x)=c_{\ga i}x^{\gg_\ga}+\ldots\qquad(c_{\ga i}\in\cc \textup{
nonzero}),
$$
and where $U(x,y)$ is a real analytic function with $U(0,0)\neq 0$. 

I see that looking at the Newton polygon alone gives me quite
detailed information on the structure of the zero set $\sxy$, as well
as of its level sets. Now I am going to proceed along the lines of the
example considered in the previous section.

I think of the $(x,y)$ plane, or rather of its positive quadrant, as
split into pieces by curves $y=x^{\gg_\ga}$ (Fig. 6). Note that the way
I number finite edges from right to left, numbers $\gg_\ga$ decrease
with $\ga$. 

\Picture{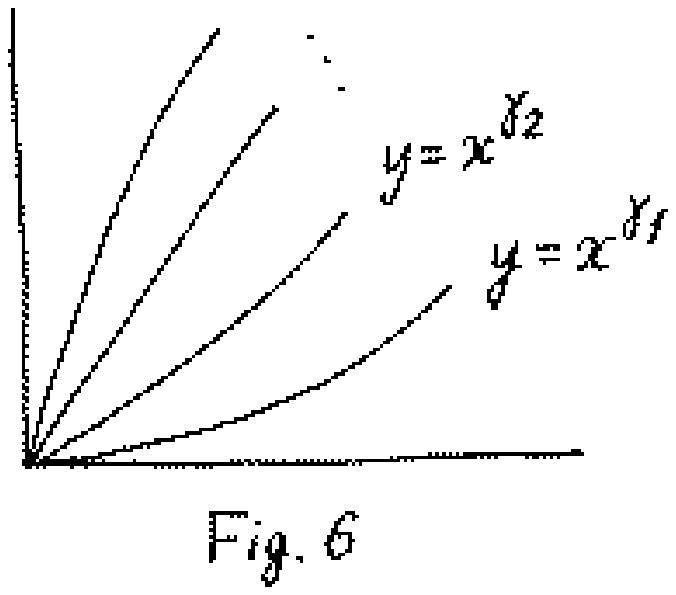}

Now I consider the dyadic partition of the positive quadrant into the
rectangles $R_{jk}$, and the corresponding smooth partition of the
$\tl$ into the operators $T_{jk}$. The rectangles $R_{jk}$ fall into two 
categories, the ones which lie far away from any of the curves
$y=x^{\gg_\ga}$, and the ones which lie close to one of these curves.

{\bf Far away rectangles. }

 Consider the rectangles lying between
$y=x^{\gg_\ga}$ and $y=x^{\gg_{\ga+1}}$. To simplify the notation, I
put $\ga=1$, but I do not assume that I am dealing with the rightmost
finite edge. The corresponding pairs of
$(j,k)$ are singled out by the condition
\beq\label{cond1}
k-\gg_1 j\<-K,\qquad k-\gg_{2}j\>K
\eeq
($K$ a large constant).

It follows from \rf{Pu} by straightforward agebra that on such a
rectangle
$$
|\sxy|\sim 2^{-jA_{2}-kB_{2}},
$$
which is a complete analogue of estimate \rf{osc11}.

Now if $B_{2}>A_2$, then I am again in the situation of
Case I of the example from the previous section. I will split the
rectangles into almost orthogonal families $k-[\gg_1j]=-r$, resum, and
get the estimate
\beq\label{gd1}
\left\|\sum T_{jk}\right\|\lesssim\gl^{-\GD_1},
\eeq
where the sum is taken over $(j,k)$ satisfying \rf{cond1}, and 
\beql{t1}
\GD_1=\frac 1{2t_1+2},
\eeq
where $(t_1,t_1)$ is the point of intersection of the straight line
passing through the edge $\ga=1$ with the bisectrix (Fig. 7). Notice
that I do not assume that $\ga=1$ is the main edge of the Newton polygon.

\Picture{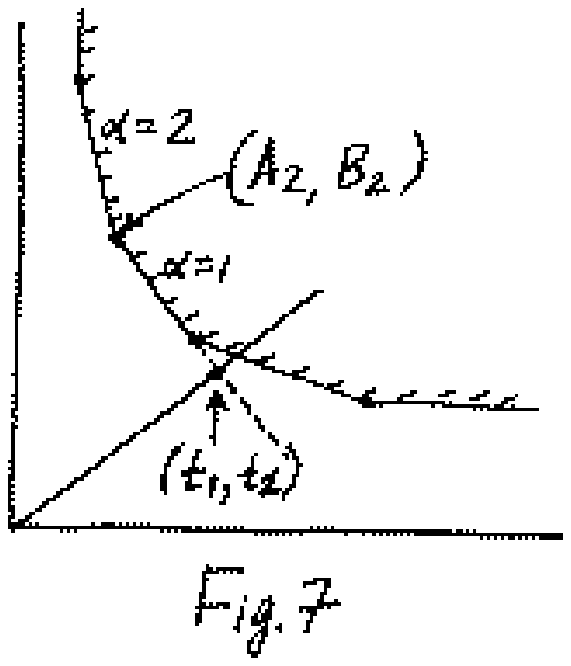}

Analogously if $B_2<A_2$, I find myself in the Case II situation. So I
will split the rectangles into the families $k-[\gg_2 j]=r$, resum, and
get the estimate \rf{gd1} with $\GD_2$ instead of $\GD_1$.

If $A_2=B_2$, I as before split into families $j+k=r$, and get the
same estimates.

Fig. 8 shows the direction of resummation for all regions of the
quadrant. Here $\gg_*$ denotes the exponent, corresponding to the main
edge (the one intersecting the bisectrix).

\Picture{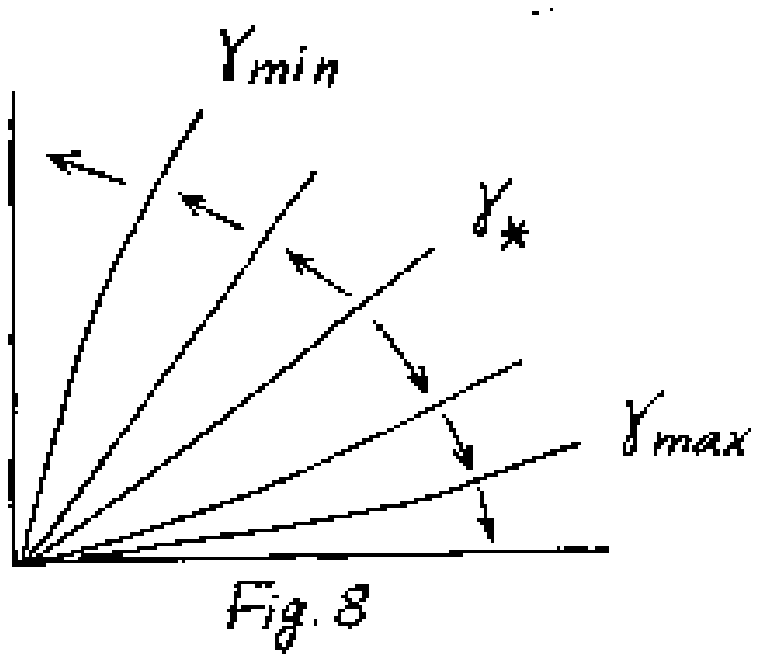}

Since obviously $\GD_\ga\>\GD$ for all edges, the above discussion
results in the needed estimate
$$
\left\|\sum_{\textup{far away }R_{jk}}T_{jk}\right\|\lesssim\gl^{-\GD}.
$$

{\bf Rectangles which are close.}

I look at the rectangles close to the curve $y=x^\gg$, where $\gg$ is
one of exponents $\gg_\ga$. These rectangles satisfy the condition
$|k-\gg j|\<K$ and form an almost orthogonal family. So it is
sufficient to prove the bound
\beq\label{bi}
\|T_{jk}\|\lesssim\gl^{-\GD}
\eeq
for each of these rectangles individually. The argument I give below 
is different from and shorter than the original proof of this
estimate given by Phong and Stein (see \cite{PS}, pp. 126--148).
 
It is easy to see from \rf{Pu} that on $R_{jk}$ the $\sxy$ has the
following behavior
\beq\label{onrjk}
\sxy\sim const. 2^{-jA_2-j\gg B_1}\prod_{i=1}^n(y-Y_i(x)),
\eeq
where to simplify the notation I put $\ga=1$, so that the edge in
question joins $(A_1,B_1)$ with $(A_2,B_2)$, but I am not going to assume
that this is the rightmost edge. I also dropped the index $\ga$ from
$n_\ga$, $\gg_\ga$, and $Y_{\ga i}(x)$.

The branches $y=Y_i(x)$ are in general complex-valued.
I introduce 
$$
Z_i(x)=\Re Y_i(x).
$$
The $Z_i(x)$ are smooth analytic functions, and for $x\sim 2^{-j}$ I
have
\beql{L}
\frac d{dx}Z_i(x)\sim x^{\gg-1}\sim2^{-j(\gg-1)}\sim const =: L.
\eeq

Some of the curves $y=Z_i(x)$ may intersect the rectangle $R=\rjk^*$
(Fig. 9). The geometry of the problem suggests to take a Whitney
decomposition of the set $R\backslash Z$, where
$$
Z=\bigcup_i \{(x,y):y=Z_i(x)\},
$$
into rectangles $R_l$ of sixe $2^{-m_l}\times L2^{-m_l}$ such that the
distance from $R_l$ to Z is $\sim 2^{-m_l}$ in anistotropic norm
$|x|+L^{-1}|y|$.

\Picture{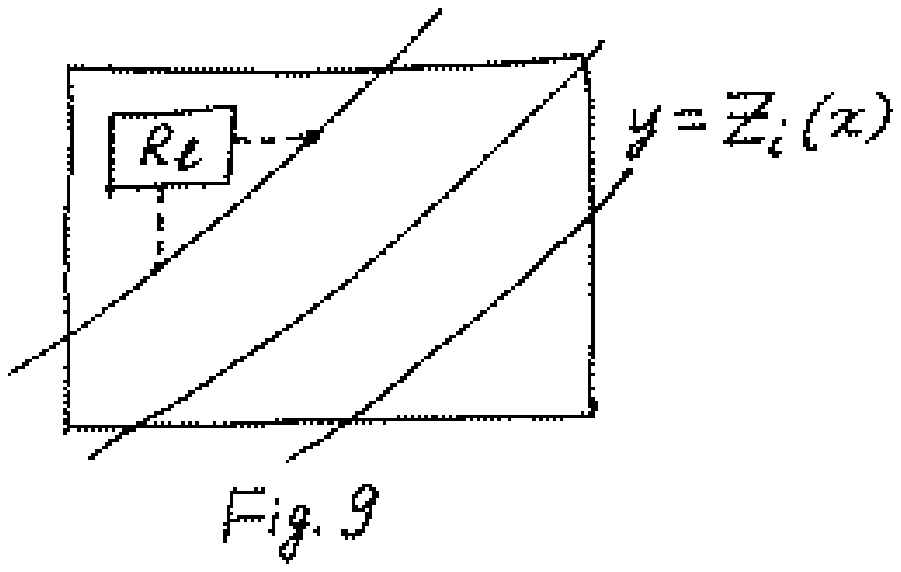}

The easiest way to arrange such a decomposition is to dilate the
picture vertically by the factor of $L^{-1}$, do the usual Whitney
decomposition, and contract back.

A quite obvious but important point follows from \rf{L}: For each $m$,
the subfamily of the rectangles $R_l$ having $m_l=m$ is almost
orthogonal.

Now I am going to smoothly localize $T_{jk}$ to $R_l^*$'s, denoting
the corresponding partial operators $T_l$. This does not break almost
orthogonality, and I have
\beql{ssum}
\|\tjk\|\<\sumi_{m=j}\left\|\sum_{l:m_l=m}T_l\right\|\lesssim\sumi_{m=j}\sup_{l:m_l=m}\|T_l\|.
\eeq   

I now turn to estimating the norm of $T_l$. Note that on $R_l^*$ I
have
$$
|y-Y_i(x)|\sim|y-Z_i(x)|+2^{-j\gb_i}\sim
 2^{-j(\gg-1)-m_l}+2^{-j\gb_i},
$$
where $x^{\gb_i}$ is the first nozero term in the expansion of $\Im
Z_i(x)$ ($\gb_i:=\infty$ if this expansion is identically zero).

It follows that on $R_l^*$
$$
|\sxy|\sim2^{-j A_2 -j\gg
B_1}\prod_{i=1}^n(2^{-j(\gg-1)-m_l}+2^{-j\gb_i})=:\mu.
$$
It is also not difficult to see that supplementary conditions of Lemma \ref{oscest}
$$
|\D^N_y\sxy|\lsim \mu 2^{j(\gg-1)N}\qquad(N=1,2)
$$
are satisfied on $R_l^*$.

Now I am prepared to apply the oscillatory estimate to $T_l$. Namely, 
Lemma \ref{oscest} gives 
\beql{osc22}
\|T_l\|\lsim
(\gl\mu)^{-1/2}\<\gl^{-1/2}2^{j(A_2+B_2\gg)/2}2^{(m_l-j)n/2}.
\eeq
(The last inequality follows by taking the lower bound for $\mu$
ignoring contributions of $\gb_i$, and also by using $n=B_2-B_1$.)

As always, I also have the following size estimate:
\beql{size22}
\|T_l\|\lsim2^{-m_l-j(\gg-1)/2}.
\eeq

The natural decay in $\gl$ that I expect is $\gl^{-\GD_1}$, where
$\GD_1$ has the same meaning as in \rf{gd1}, in particular,
$$
\GD_1=\frac{1+\gg}{2(A_2+1)+2(B_2+1)\gg}.
$$
So I take the geometric mean of \rf{osc22} and \rf{size22} with the
corresponding exponents $\theta=2\GD_1$ and $1-\theta$. As the reader
may check, I get (see \rf{t1} for the definition of $t_1$)
\beql{est22}
\|T_l\|\lsim\gl^{-\GD_1} 2^{-\tilde{m}(2t_1-n)\GD_1},\qquad \tilde
m=m_l-j\>0.
\eeq 

Now please note that I may assume $\gg\>1$. Indeed, if $\gg<1$, then I
just switch to the adjoint of $\tjk$, which amounts to interchanging
roles of $x$ and $y$ and transforms $\gg\to 1/\gg$. 

Further, if $\gg>1$ strictly, then it is easy to see geometrically
that necessarily $t_1>n/2$ no matter how the edge of the Newton
polygon lies. In this situation I can substitute \rf{est22} into
\rf{ssum}, sum in $\tilde m$, and get the required estimate \rf{bi}
(note that $\GD_1\>\GD$).

The special case $t_1=n/2$ can happen only if $\gg=1$ and the Newton
polygon has only one finite edge joining points $(n,0)$ and $(0,n)$.
In this case I avoid taking the geometric mean and substitute
\rf{osc22} and \rf{size22} directly into \rf{ssum}:
\beql{sum22}
\|\tjk\|\lsim \sumi_{\tilde m=0}\min(\gl^{-1/2}2^{jn/2+\tilde m
n/2},2^{-\tilde m-j})\lsim \gl^{-1/(n+2)},
\eeq
as a simple analysis shows. (Find $\tilde m_*=\tilde m$ for which the
progressions balance. Consider the cases $\tilde m_*<0$ and $\tilde
m_*\>0$.) 
This is the right estimate in this
particular case. \qedsymbol

\subsection {Discussion and outlook}

This finishes the proof of the upper bound in the real analytic case.

The main components of the proof, such as 
\begin{itemize}
\item the use of the Puiseux expansion,
\item the resummation procedures used to estimate far away from the branches,
\item the Whitney decomposition method used close to the branches
\end{itemize}
are going to carry over to the $\cin$ case either verbatim or with
small modifications, as the reader will see in the coming chapters.

Jumping slightly ahead of time, I am going to say that the only
crucial difference, actually the one responsible for the presence of
$\log\gl$'s in Theorem 1.1, is going to come from the possible occurence of
multiple real nondifferentiable branches. Namely, in the $\cin$ case I
may have a situation like the one shown in Fig. 10, when the zero set
of $\sxy$ has several branches, which, although being close to each
other to infinitely high order, are nevertheless non-coinciding and in
fact nondifferentiable. This of course would be impossible in the real
analytic case. 

\Picture{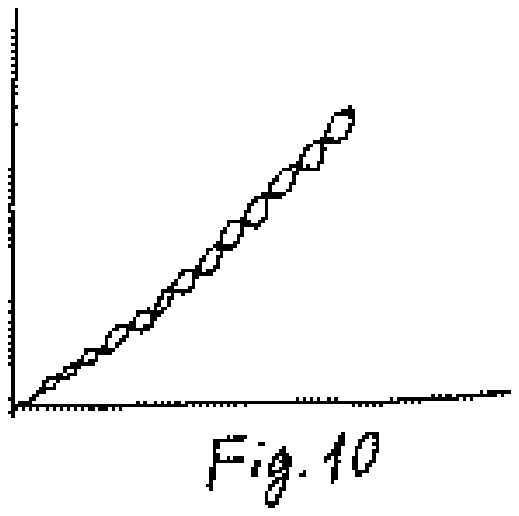}

The problem with such a situation is that I lose condition \rf{L},
which was the basis of almost orthogonality, and
the whole Whitney decomposition procedure is going to become useless
near these multiple branches.

The way I am going to fight this difficulty will be to localize away
from the branches by a very narrow cutoff. The Whitney
decomposition will still work away from the branches, and near the
branches I will have to use a completely different argument, based on a method
due to Seeger \cite{See1}.

\chapter{Smooth Puiseux theorem}

The next 2 chapters are devoted to the proof of Theorem \ref{main}.
While in the previous chapter, which was supposed to be expository,
I was allowing myself to be informal at times, from now on I will
strive to provide full details.
    
In this chapter I explain how the Puiseux expansion \rf{Pu}
generalizes from the real analytic to the $\cin$ case.

\section{Algebraic notation}

Proving theorems about $\cin$ functions often involves an intermediate
step, when the analysis is done purely algebraically within the
category
of formal power series. I am going to employ this very
strategy. Here, I will set up some algebraic notation.

First recall that for any ring $R$, the symbols $R[t]$ and $R[[t]]$ denote the rings of 
polynomials and, respectively, formal power series in indeterminate $t$ with coefficients from $R$. 
This notation can be iterated, e.g. $R[[x]][y]$ is the ring of
polynomials
in $y$ with coefficients which are elements of $R[[x]]$, and
$R[[x,y]]$ is the ring of double formal power series. 

Factorization formulas for $\cin$ function, which I am going to prove
in this chapter, are going to be valid in a small
neighborhood of the origin. Since I do not care how small this
neighborhood is, it will be convenient to formulate the results for
\emph{function-germs} rather than functions. 

An identity involving
several function-germs is defined to be true if there exist
functions from the equivalence classes of these germs such that in the
intersection of their domains of definition the identity is true in
the usual sense.

Basically, this convention will
spare me the necessity to repeat the phrase ``There exists a small
neighborhood of the origin $U$ such that in $U$ \ldots'' every time.

I will make use of the following rings of germs of $\cc$-valued
functions:
\begin{itemize}  
\item $C\(x\)$ --- continuous functions at the origin of $\rr$; 
\item $\cin\(x\)$ and $\cinxy$ --- $\cin$ functions at the origin of $\rr$
and $\rr^2$, respectively;
\item
 $C_+\(x\)$ and $\cin_+\(x\)$ --- rings of one-sided germs; consist of
(the equivalence classes of) functions $f(x)$ 
defined in a left half-neighborhood of zero of
the form $[0,\eps)$, where $\eps>0$ can depend on $f(x)$, which are
continuous, respectively $\cin$, up to zero;
\item 
$A_+\(x^\gg\)$, $\gg>0$, --- the subring of 
$C_+\(x\)$ consisting of germs $f(x)$, for which there exists a
series
$\bar f (x)\in\cc[[x^\gg]]$, $\bar f(x)=\sumi_{n=0}c_n x^{n\gg}$,
such that $f(x)\sim\bar f(x)$ in the sense that for any $N$
\[
f(x)-\sum_{n=0}^N c_n x^{n\gg}=O(x^{(N+1)\gg}),\qquad x\to0.
\]
Such an $\bar f(x)$ is uniquely determined and is called the
\emph{asymptotic expansion} of $f(x)$.
\end{itemize}

Notice that for the elements of $\cinxy$, $\cin\(x\)$, and
$\cin_+\(x\)$, I can talk about their Taylor series at the origin.
A germ whose Taylor series is zero is called \emph{flat}.

The rings of germs of $\rr$-valued functions will be denoted by adding
an $\rr$ to the above notation, e.g. $\rr\cinxy$.

\section{Puiseux decomposition of $\cin$ functions}

Now I am going to state the main result of this chapter. The proof
will be given in the following sections.

Let $F(x,y)\in\rr\cinxy$. Denote by $\GG=\GG(F)$ the Newton polygon of $F(x,y)$,
and assume that $\GG\ne\emptyset$, so that $F$ is not flat. 

Let $\ga$ run through all compact edges of the boundary of $\GG$. For each edge $\ga$ joining
integer points $(A_\ga,B_\ga)$ and $(A_\ga',B_\ga')$, where $B_\ga'>B_\ga$,
put 
$$
n_\ga=B_\ga'-B_\ga,\qquad \gg_\ga=\frac{A_\ga-A_\ga'}{B_\ga'-B_\ga}.
$$
Let also $A$ be the $x$-coordinate of the vertical infinite edge, 
and $B$ be the $y$-coordinate of the horizontal infinite edge
of $\GG$. 

\begin{Prop}\label{31} In the above conditions, the germ $F(x,y)$ admits in the region $x,y>0$ a
factorization of the form
\begin{equation}
\label{factor1}
F(x,y)=U(x,y)\prod_{i=1}^{A}(x-X_i(y)) 
\prod_{i=1}^{B}(y-Y_i(x)) 
\prod_\ga\prod_{i=1}^{n_\ga}(y-Y_{\ga i}(x)),
\end{equation}
where\\
\textup{(1)} $U(x,y)\in\rr\cinxy$, $U(0,0)\ne 0$,\\   
\textup{(2)} all $X_i(x), Y_i(x)\in C_+\(x\)$, and
$X_i(x), Y_i(x)=O(x^N)$ as $x\to 0$ for any $N>0$,\\ 
\textup{(3)} all $Y_{\ga i}(x)\in A_+\(x^{1/n!}\)$ for $n=
B+\sum_\ga n_\ga$ with asymptotic expansions of the form $Y_{\ga i}(x)=
c_{\ga i} x^{\gg_\ga}+\ldots$ as $x\to 0$, where $c_{\ga i}\ne 0$,\\
\textup{(4)} if $Y(x)$ is any of the functions $Y_{\ga i}(x)$, and if 
$f(x,y)=\prod (y-Y_{\ga i}(x))$ is the product over all $i$ such that 
$Y_{\ga i}(x)$ has exactly the same asymptotic expansion as $Y(x)$,
then $f(x^{n!},y)\in \cin_+\(x\)[y]$,\\
\textup{(5)} if in \textup{(4)} I additionally assume that the
asymptotic expansion of $Y(x)$ is real, then $f(x,y)$ is also real.
\end{Prop}

\vspace{10pt}
This result copies \rf{Pu} in the part that concerns the number of
branches and the leading terms in their asymptotic expansion. However,
there are also substantial differences, such as:
\begi
\item in general, branches $X_i(y)$ and $Y_i(x)$ infinitely tangent to coordinate
axes but not coinciding with them are present;
\item I cannot claim that branches $Y_{\ga i}(x)$ are
differentiable; (4) is the best that is true in general.
\ei
These differences are for real, as very simple example show. For
instance, one can take $\sxy=y^2+a(x)y+b(x)$ with a choice of
coefficients so that the determinant oscillates around zero as $x\to 0$.

To the best of my knowledge, Puiseux decompositions of $\cin$
functions in the form of Proposition 3.1 or of a similar kind have not 
appeared in the literature before. However, granted Lemma 3.2 below,
the proof of Proposition 3.1 follows a rather standard path, well known
say in the singularity theory of $\cin$ and analytic functions, see
e.g. Arnold et.al.~\cite{Arnold}, or Artin \cite{Artin}.  
  
\section{Preparation to the proof}
 
The proof relies on the following result, which is well known in
the theory of plane algebraic curves under the same generic name of
the Puiseux theorem. A proof can be found in \cite{Wal}, p.~98ff, or \cite{Bour}, A.V.150.

\begin{Lem}
\label{puiseux} Let $\bar F(x,y)\in\cc[[x]][y]$ be of the form
\[
\bar F(x,y)=y^n+
\bar c_{n-1}(x) y^{n-1}+\ldots+\bar c_0(x),\qquad \bar c_i(x)\in\cc[[x]],
\]
where the zeroth order terms of all $c_i(x)$ vanish. Let $\ga$, $n_\ga$, $\gg_\ga$, $B$
be
defined via the Newton polygon $\GG=\GG(\bar F)$ in the same way as in the
Proposition. Then there exists a factorization
\begin{equation}
\label{puiseux1}
\bar F(x,y)=y^{B}\prod_\ga\prod_{i=1}^{n_\ga}(y-\bar Y_{\ga i}(x)),
\end{equation}
where the series $\bar Y_{\ga i}(x)\in\cc[[x^{1/n!}]]$ are of the
form $\bar Y_{\ga i}(x)=c_{\ga i} x^{\gg_\ga}+\ldots$ with $c_{\ga
i}\ne 0$.
\end{Lem}

The following lemma will be used to pass from factorizations of
formal power series (obtained via Lemma \ref{puiseux}) to factorizations of
function-germs in the $\cin$ category. The proof uses standard
technology usually applied in such situations.

\begin{Lem} 
\label{poly}
Let $P(x,y)\in\cin\(x\)[y]$ be of the form
\[
P(x,y)=y^n+
 c_{n-1}(x) y^{n-1}+\ldots+c_0(x),\qquad c_i(x)\in\cin\(x\),
\]
where $c_i(0)=0$ for all $i$. Let $\bar P(x,y)\in\cc[[x]][y]$ 
be the formal Taylor series of $P(x,y)$ at the origin. 
Let $\bar Y(x)\in\cc[[x]]$ be a root of multiplicity $m$, $1\<m\<n$,
of $\bar P(x,y)$ considered as a polynomial in $y$, which means that
\[
\bar P(x,\bar Y(x))=\ldots=\bar P_y^{(m-1)}(x, \bar Y(x))=0,
\qquad \bar P_y^{(m)}(x, \bar Y(x))\ne 0
\]
as elements of $\cc[[x]]$. Then there exist $m$ function-germs 
$Y_1(x),\ldots,Y_m(x)\in C\(x\)$ such that\\
\textup{(1)} all $Y_i(x)\sim \bar Y(x)$ as $x\to 0$,\\
\textup{(2)} all $P(x,Y_i(x))=0$ for $x>0$,\\
\textup{(3)} $\prod_{i=1}^m(y-Y_i(x))\in\cin\(x\)[y]$,\\
\textup{(4)} if we additionally assume that $P(x,y)$ and $\bar Y(x)$
are real, then $\prod_{i=1}^m(y-Y_i(x))$ is also real.
\end{Lem}
\begin{proof} Let $\tilde Y(x)$ be a $\cin$ function with the formal
Taylor series $\bar Y(x)$, supplied by E.~Borel's theorem. Denote 
\[
\delta_i(x)=\frac 1{i!}P_y^{(i)}(x,\tilde Y(x)),\qquad i=0,\ldots,n.
\]
Let the (nonzero by assumption) series $\bar P^{(m)}_y(x,\bar Y(x))$
starts with a term $cx^s$, $c\ne 0$, $s\in \zp$. 
Then I have $\gd_m(x)=cx^s+o(x^s)$ as $x\to 0$. On the other
hand, the functions $\gd_0(x),\ldots,\gd_{m-1}(x)$ are flat.

I will be looking for $Y_i(x)$ of the form 
\[
Y(x)=\tilde Y(x)+\delta_m(x)\ga(x),
\]
where $\ga(x)$ is an unknown continuous $\cc$-valued function-germ such that
$\ga(x)=O(x^N)$ as $x\to0$ for any $N>0$.

By Taylor's formula, the equation $P(x,Y(x))=0$ can be written as
\begin{equation}
\sum_{i=0}^n \gd_i(x)[\gd_m(x)\ga(x)]^i=0.
\end{equation}
For small $x$, this is equivalent to the equation
$w(x,\ga(x))=0$ for the function $w(x,z)$ given by
\[
w(x,z)=\sum_{i=0}^n \gd_i(x)[\gd_m(x)]^{i-m-1} z^i.
\]
Note that if $f,g\in\cin$, and $f$ is flat at the origin, while $g$ is
not flat, then $f/g$ is $\cin$ near the origin and is flat. So the
performed division by $[\gd_m(x)]^{m+1}$ is legitimate, and 
$w(x,z)\in\cin\(x\)[z]$. 

 On the complex circle $|z|=x^N$,
$N>0$, the term $z^m$ will dominate the other terms in $w(x,z)$ if $x$ is
sufficiently small. By Rouche's theorem it follows that the equation
$w(x,z)=0$ has for small fixed $x$ exactly $m$ roots in the disc
$|z|<x^N$, which I denote
$\ga_i(x)$, $i=1,\ldots,m$. I can arrange so that $\ga_i(x)$ are
continuous in $x$, and the previous argument shows that
$\ga_i(x)=O(x^N)$ for any $N>0$.

We now prove (3). Since the functions $\ga_i(x)$ 
enter the product
in (3) in a symmetric way, it is sufficient to prove that the
elementary symmetric polynomials $s_1,\ldots, s_m$ in $\ga_i(x)$ are in
$\cin\(x\)$. By the Newton relations (see
\cite{Bour}, A.IV.70), it is sufficient to prove the same for the
functions
\[
p_k(x)=\sum_{i=1}^m [\ga_i(x)]^k,\qquad k=1,\ldots,m.
\]
However, by Cauchy's formula I have that for small $x$
\[
p_k(x)=\frac1{2\pi i}\oint_{|z|=\eps} \frac {z^k w'_z(x,z)}{w(x,z)}\,dz,
\]
from where it is clear that $p_k(x)\in\cin\(x\)$, since nothing
dramatic happens to $w(x,z)$ on the circle $|z|=\eps$.

To prove (4), I notice that under the additional assumption made I
can take $\tilde Y(x)$ to be real. Then $w(x,z)\in\rr\cin\(x\)[z]$,
and therefore non-real roots $\ga_i(x)$ will appear in conjugate
pairs. Then all $p_k(x)$ will be real, which implies (4). 
\end{proof}

Now I am going to combine two previous lemmas to prove

\begin{Lem} 
\label{pfactor}
Proposition \ref{31} is true if $F(x,y)\in\rr\cin\(x\)[y]$.
\end{Lem} 
\begin{proof}
By Lemma \ref{puiseux}, the Taylor series $\bar F(x,y)\in\rr[[x]][y]$ 
of $F(x,y)$ has a factorization \rf{puiseux1}. Consider the function 
$P(x,y)=F(x^{n!},y)$. Its Taylor series 
has the form $\bar P(x,y)=\bar F(x^{n!},y)$, and so factorizes as
\[
 \bar P(x,y)=y^{B}\prod_\ga\prod_{i=1}^{n_\ga}(y-\bar Y_{\ga i}(x^{n!})).
\]
Let $\bar Y(x)$ be one of the series $\bar Y_{\ga i}(x^{n!})\in\cc[[x]]$,
and assume that among all the $\bar Y_{\ga i}(x^{n!})$ there are exactly $m$
series coinciding with $\bar Y(x)$. Then $y=\bar Y(x)$ is a root of
multiplicity $m$ of the polynomial $\bar P(x,y)\in\rr[[x]][y]$, and by
Lemma \ref{poly} I conclude that there exist $m$ functions $Y_i(x)\in
C\(x\)$, $i=1,\ldots,m$, such that (1)--(3) from the formulation of
the lemma are true. 

In view of (3), we can divide $P(x,y)$ by
$\prod(y-Y_i(x))$, and the result is again a polynomial $\tilde
P(x,y)$ from $\cin\(x\)[y]$. The Taylor polynomial of $\tilde
P(x,y)$ will be $\bar P(x,y)$ divided by $(y-\bar Y(x))^m$. Now we can
apply Lemma \ref{poly} to $\tilde P(x,y)$ choosing a different $\bar
Y(x)$ etc. 

By repeating this operation several times, I get a complete factorization of $P(x,y)$. The
required factorization of $F(x,y)$ is then obtained by the inverse
substitution $x\mapsto x^{1/n!}$. The property (5) is ensured by 
splitting off all real series $\bar Y(x)$ before non-real ones in the
above argument.  
\end{proof}

Proposition \ref{31} will be reduced to Lemma \ref{pfactor}
by means of the following \emph{Malgrange preparation theorem} (see
\cite{GG}, p.~95).

\begin{Lem}
\label{malgrange}
Let $F(x,y)\in\rr\cinxy$, and assume that $F(0,y)$ is not flat, so that
$F(0,y)=cy^n+o(y^n)$, $y\to0$, for some $n\in\zp$ and $c\ne0$. Then
there is a factorization 
\[
F(x,y)=U(x,y)P(x,y),
\]
where \\
\textup{(1)} $U(x,y)\in\rr\cinxy$, $U(0,0)\ne 0$,\\
\textup{(2)} $P(x,y)\in\rr\cin\(x\)[y]$ is of the form
\[
P(x,y)=y^n+c_{n-1}(x)y^{n-1}+\ldots+c_0(x),
\]
where all $c_i(x)\in\rr\cin\(x\)$, $c_i(0)=0$.
\end{Lem}

\section{Proof of the Proposition}

Notice that the Newton polygon is invariant with \mbox{respect} to
multiplication by a nonzero $\cin$ function (see Phong and Stein
\cite{PS}, p.~112). Therefore, for the functions $F(x,y)$ such that $F(0,y)$
is not flat (which is equivalent to having $A=0$) the proposition
follows immediately from Lemmas \ref{malgrange} and \ref{pfactor}.

Assume now that $A>0$. In this case we must somehow separate the roots
infinitely tangent to the $y$-axis. This can be done as follows. Since
$F(x,y)$ is not flat at the origin, there exists a rotated orthogonal
system of coordinates $(x',y')$ such that the restriction of $F$ to
the $y'$-axis is not flat. So we can apply Lemma \ref{malgrange} to
$F$ written in coordinates $(x',y')$. Let $P(x',y')$ be the arizing
polynomial.    

If $y'=ax'$ is the equation of the old $y$-axis in the new
coordinates, then $y'=ax'$ will be a root of multiplicity $A$ of $\bar
P(x',y')\in\rr[[x']][y']$. So we can apply Lemma \ref{poly} and obtain
$A$ roots $y'=Y_i(x')$, $i=1,\ldots,A,$ of $P(x',y')=0$, such that
$Y_i(x')\sim ax'$.

Moreover, by Lemma \ref{poly} (3),(4) we will have that
$Q(x',y')=\prod(y'-Y_i(x'))$ is in $\rr\cin\(x'\)[y']$. So we can
divide $P(x',y')$ by $Q(x',y')$, and the quotient will be a $\cin$
function, which is no longer flat on the old $y$-axis.

Let $\tilde F(x,y)$ be this last quotient written in the old system of
coordinates. Then the Newton polygon of $\tilde F$ is just $\GG(F)$ shifted $A$
units to the left. So we can factorize $\tilde F(x,y)$ as in the case
$A=0$ described above. 

It remains to get a factorization of $Q(x',y')$ in the old
coordinates. It is clear that the Taylor series of $Q$ written in
the coordinates $(x,y)$ consists of one term $cx^A$. Interchanging
the roles of $x$ and $y$ brings us back to the case $A=0$, and the
required factorization of the form $\prod(x-X_i(y))$ can be obtained
as described above.

\chapter{Upper bound for smooth case}

In this chapter I am going to prove Theorem \ref{main}.

\section{Beginning of the proof}

The proof starts just like in the real analytic case.

I decompose the operator $\tl$ as
\[
\tl=\sum_{\pm}\sum_{j,k} T_{jk}^{\pm\pm},
\] 
where $T_{jk}^{++}$ is defined as 
\[
T_{jk}^{++}f(x)
=\int_{-\infty}^{\infty}e^{i\gl S(x,y)}\chi_j(x)\chi_k(y)\chi(x,y)f(y)\,dy. 
\] 
Here $\sum_j\chi_j(x)=1$ is a smooth dyadic partition of unity on
$\rr^+$,
so that the kernel of $T_{jk}^{++}$ is supported on the rectangle
$R_{jk}=[2^{-j-1},2^{-j+1}]\times[2^{-k-1},2^{-k+1}]$. Three other
$\pm$
combinations refer to the quadrants defined by specific signs of $x$
and $y$. We restrict ourselves with the positive quadrant, the other
ones
being exactly similar, and denote $T_{jk}^{++}$ by simply $T_{jk}$.
 
By Proposition \ref{31} applied to $F(x,y)=\sxy(x,y)$, there is 
a neighborhood of the origin $V$ such that in $V\cap \rrp^2$ 
there exists a factorization of the form \rf{factor1}. 
I assume that $\supp \chi\subset V$. The singular variety 
$$
\cal Z=\{(x,y)\in V:F(x,y)=0\}
$$ 
now splits into branches corresponding to the factors in the RHS of \rf{factor1}. Note, however, that some
of these branches may contain an imaginary component.

Let $R_{jk}^*$ denote the double of $R_{jk}$.
I fix a large constant $D$ such that if the pair $(j,k)$
satisfies the condition $\min_\ga|k-j\gg_\ga|\>D$, then $y-c_{\ga
i}x^{\gg_\ga}\ne 0$ on $R_{jk}^{**}$ for all $c_{\ga i}x^{\gg_\ga}$
occurring
as the lowest order terms of the asymptotic expansions of $Y_{\ga i}(x)$
in Proposition \ref{31}. 

Let me number the compact edges $\ga$ of the boundary of the Newton polygon $\GG(F)$ 
from right to left, so that 
$$
\gg_1>\gg_2>\ldots>\gg_{\ga_0},
$$ 
where $\ga_0$ is the total number of compact edges. Also put $\gg'=\gg_{\ga_0}/2$ if
$A>0$, $\gg'=0$ otherwise; $\gg''=2\gg_1$ if
$B>0$, $\gg''=\infty$ otherwise. 

Consider the following splitting of $\tl$:
\beql{split}
\tl= \left(T'+T''+\sum _{\nu=1}^{\ga_0-1} T_\nu\right)+ (T_x +T_y) +
\sum_{\nu=1}^{\ga_0} T^\nu.
\eeq
Here  
\beql{talpha}
T_\nu=\sum_{j\gg_{\nu+1}\ll k\ll
j\gg_{\nu}} T_{jk},\qquad  1\<\nu\<\ga_0-1,
\eeq
\[
T'=\sum_{j\gg'\ll k\ll
j\gg_{\ga_0}} T_{jk},\qquad
T''=\sum_{j\gg_{1}\ll k\ll
j\gg''} T_{jk},
\]
($a\ll b$ stands for $a\<b-D$) constitute the part of $\tl$ supported relatively
far away from $\cal Z$. Further,
\[
T_y=\sum_{k\ll\gg'j} T_{jk},\qquad T_x=\sum_{\gg''j\ll k} T_{jk}.
\]
constitute the part of $\tl$ supported near the branches of
$\cal Z$ which are infinitely tangent to the coordinate axes. Finally
\begin{equation}
\label{tnu}
T^\nu=\sum_{\gg_\nu j-D<k<\gg_\nu j+D} T_{jk},\qquad  1\<\nu\<\ga_0,
\end{equation}
are the part of $\tl$ supported near all other branches of $\cal Z$.

In the following sections I will prove the upper bound claimed in
Theorem 1.1 for the norms of all operators in the RHS of \rf{split}: $T',T''$ and $T_\nu$
(Section 4.2), $T_x$ and $T_y$ (Section 4.3), and finally $T^\nu$
(Sections 4.4 and 4.5). This will prove Theorem 1.1. 


\section{Estimates far away from $\cal Z$}

In this section, I prove that $\|T_\nu\|\<C\gl^{-\GD}$ for each $\nu$. 
The reader will believe me that with minor modifications the
argument given below will also produce the same estimate for $T'$,
$T''$.

The proof is very similar to the argument given in Section 2.2.3 for
far away rectangles. I will just provide some extra details about
estimating the size of $F=\sxy$ on the support of $\tjk$ and about
checking conditions of Lemma 2.3. In what concerns subsequent
resummation of the individual $\|\tjk\|$ estimates, the
argument goes through verbatim.

Take an operator $T_{jk}$ entering
the RHS of \rf{talpha}. I may reduce $V$ if necessary so that on the 
part of $\cal Z$ inside $V$ the functions $Y_{\ga i}$ do not differ much
from the first terms of their asymptotic expansions. Assume that 
$T_{jk}$ is nonzero, which means that $R_{jk}\cap V\ne\emptyset$.
Then it is clear from the definition of the constant $D$ that the
factors in the RHS of \rf{factor1} can estimated as follows for 
$(x,y)\in R_{jk}$ (see Fig. 10a):
$$
|x-X_i(y)| \approx 2^{-j},
$$
\begin{equation}
\label{boundsaway}
|y-Y_i(x)| \approx 2^{-k},
\eeq
$$
|y-Y_{\ga i}(x)| \approx 
\cases{
2^{-k}, &$\ga\<\nu$,\cr
2^{-j\gg_\ga},&$\ga>\nu$.\cr}
$$
($a\approx b$ means $C^{-1}b\<a\<Cb$, where $C>0$ is an unimportant
constant independent of $j,k,\gl$).

\Picture{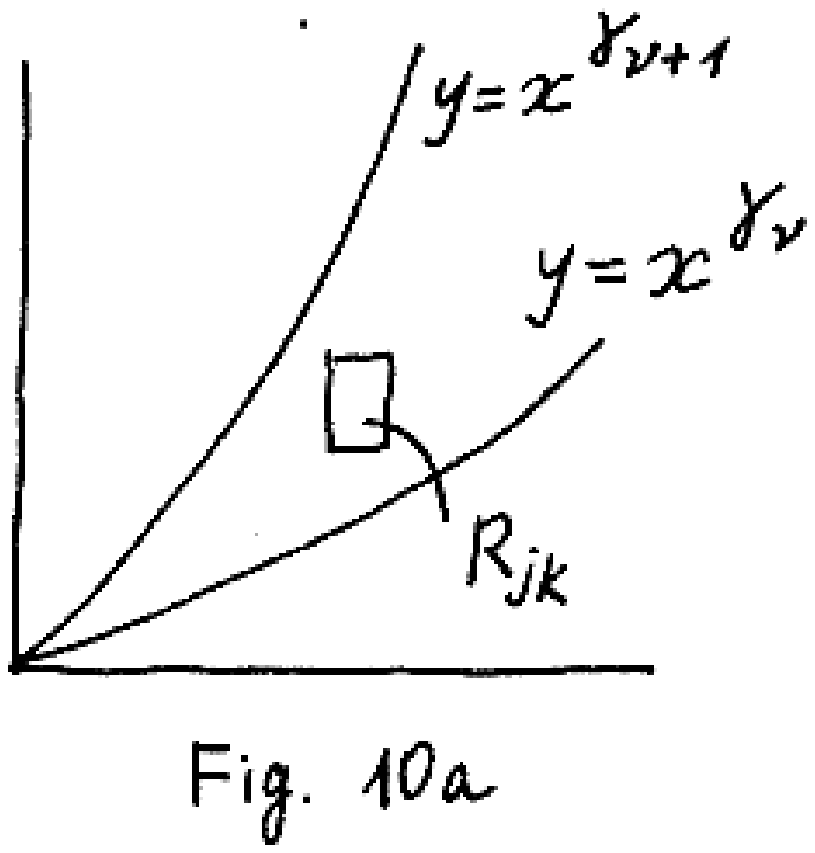}

Therefore it follows from \rf{factor1} that on $R_{jk}$
\begin{equation}
\label{hessian}
|F|\approx 2^{-jA}2^{-kB}\prod_{\ga\<\nu} 2^{-kn_\ga}
\prod_{\ga>\nu}
2^{-j\gg_\ga n_\ga}=:\mu.
\end{equation}
The numbers $\gg_\ga$, $n_\ga$ can be found from the Newton polygon 
$\GG(F)$ as described in Proposition 3.1. Using this 
information, I find that
\[
\mu=2^{-jA'_\nu-kB'_\nu}.
\]

I further claim that on $R_{jk}$
\begin{equation}
\label{derivatives}
|\D_y^n F|\lesssim \mu2^{kn},\qquad n=1,2.
\end{equation}
Indeed, when differentiating the RHS of \rf{factor1} in $y$, the derivative
can fall one either $U(x,y)$, or $\prod(x-X_i(y))$, 
or one of the 
remaining terms. In the first case, I simply get a bounded factor. In the
second
case, I get an even better factor of $O(2^{-kN})$ for any $N>0$, since 
the product in question is a $\cin$ function whose Taylor series
at the origin is $x^{A}$. Finally, in the third case I get a
factor of the form $(y-Y_i(x))^{-1}$ or $(y-Y_{\ga i}(x))^{-1}$, which
is $O(2^k)$ in view of \rf{boundsaway}. This argument works equally
well
for the second derivative, giving \rf{derivatives}.  
 
The rectangle $R_{jk}$ is of size $\gd_x\times\gd_y$ with $\gd_x\approx
2^{-j}$,
$\gd_y\approx2^{-k}$. So the conditions of Lemma \ref{oscest}
are satisfied, and I obtain the oscillatory estimate 
\begin{equation}
\label{oscillatory}
\|T_{jk}\|\lesssim \gl^{-1/2}2^{(jA'_\nu+kB'_\nu)/2}.
\end{equation}
On the other hand, the size estimate following from Lemma \ref{size} is
\begin{equation}
\label{sizeest}
\|T_{jk}\|\lesssim 2^{(j+k)/2}.
\end{equation}

The rest of the proof goes through exactly as described in Section
2.2.3. Namely, I am going to split the operators $\tjk$ constituting
$T_\nu$ 
into almost orthogonal families $k-[\gg_\nu j]=-r$, or $k-[\gg_{\nu+1}
j]=r$, or $k+j=r$, depending if $B'_\nu$ is larger, smaller, or equal
to $A'_\nu$. Then I am going to resum and get the $\gl^{-\GD}$ estimate
for $\|T_\nu\|$. I will not repeat the details.


\section{Estimates near the coordinate axes}

In this section, I will prove the estimate
$\|T_x\|\lesssim
\gl^{-\GD}$. The same estimate will be true for $T_y$, since 
taking the adjoint of $T$ brings $T_y$ to the form of $T_x$. 
I may of course assume $B\>1$, since otherwise $\gg''=\infty$ and $T_x=0$.

Notice that in the real analytic case there was no need to introduce
this special localization along the coordinate axis. In the notation
of Section 4.1, operator $T_y$ could be included into the $T'$ part
and treated along the same lines as the $T_\nu$. Analogously $T_x$
could be united with $T''$. However, in the $\cin$ case the possible
presence of the branches infinitely tangent to coordinate axes asks
for this additional localization.

I represent $T_x$ as (see Fig. 10b)

\Picture{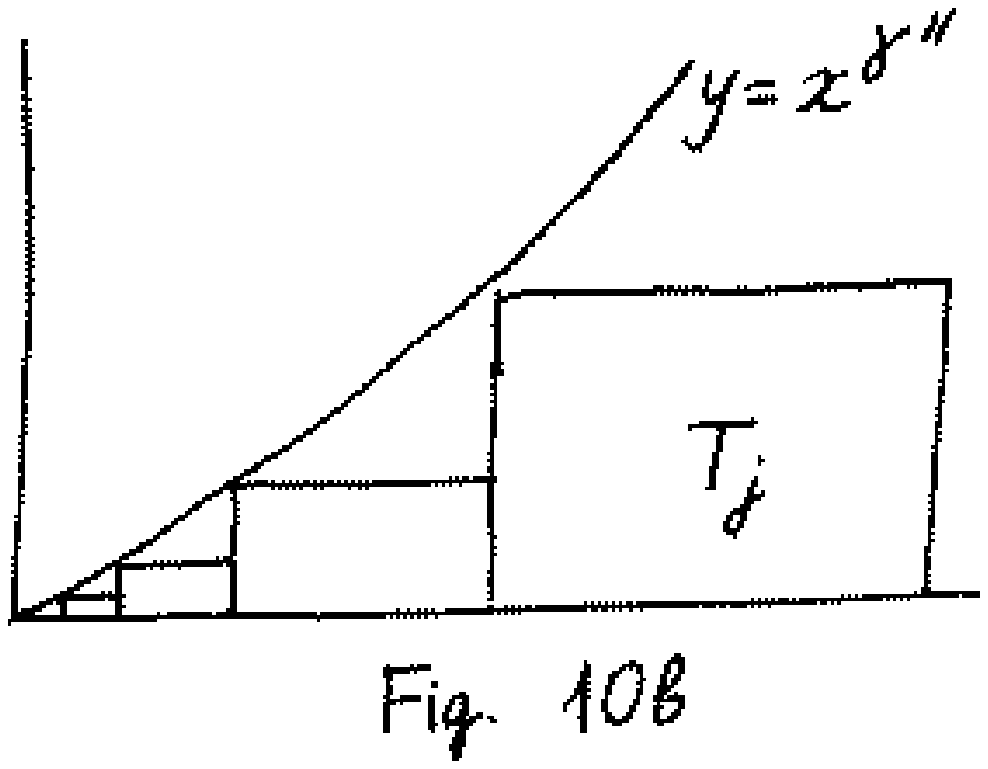}

\[
T_x=\sum_j T_j,\qquad T_j=\sum_{k:\gg''j\ll k} T_{jk},
\]
and claim that\\
(1) $\|T_j\|\lesssim\gl^{-\GD},$\\
(2) $\|T_j^*T_{j'}\|=0$ for $|j-j'|\>2$,\\
(3) $\|T_j T_{j'}^*\|\lesssim\gl^{-2\GD}2^{-\eps|j-j'|}$ for some
$\eps>0$.\\
If I prove all these, the estimate $\|T_x\|\lesssim
\gl^{-\GD}$ will follow from the Cotlar--Stein lemma.

I have 
\[
T_{j}f(x)
=\int_{-\infty}^{\infty}e^{i\gl
S(x,y)}\chi_j(x)\tilde\chi_j(y)\chi(x,y)f(y)\,dy,  
\]
where $\tilde \chi_j=\sum_{\gg''j\ll k}\chi_k$, so that the support
of $\tilde\chi$ is contained in $[0,C2^{-\gg'' j}]$. The property (2) is
obvious. Further, the operator $T_jT_{j'}^*$ has the kernel
\[
K(x_1,x_2)=\chi_j(x_1)\chi_{j'}(x_2)\int e^{i\gl[S(x_1,y)-S(x_2,y)]}
\tilde\chi_j(y)\tilde\chi_{j'}(y)\chi(x_1,y)\chi(x_2,y)\,dy.
\] 
I want to estimate this by the following variant of the standard van der
Corput lemma (see \cite{STEIN}, Corollary on p.~334).

\begin{Lem} 
\label{vandercorput}
Let $k$ be a positive integer, $\GF\in C^{k}[a,b]$, $\GP\in C^1[a,b]$, and assume that
$\GF^{(k)}\>\mu>0$ on $[a,b]$. If $k=1$, assume
additionally that $\GF'$ is monotonic on $[a,b]$. Then
\[
\biggl|\int_a^b e^{i\gl\GF(y)}\GP(y)\,dy\biggr|\lesssim
(\gl\mu)^{-1/k}\biggl(|\GP(b)|+\int_a^b |\GP'|\biggr).
\]
\end{Lem}

Assume that $j'\>j$. I apply this lemma with $[a,b]=[0,C2^{-\gg''j}]$, 
 $k=B+1\>2$,
$$
\GF(y)=S(x_1,y)-S(x_2,y),
$$
$$
\GP(y)=\tilde\chi_j(y)\tilde\chi_{j'}(y)\chi(x_1,y)\chi(x_2,y).
$$
It is clear that $|\GP(b)|+\int_a^b |\GP'|\lesssim 1$. Further (recall
that I denoted $\sxy=F$), 
\[
\GF^{(B+1)}(y)=\D_y^{B+1} S(x_1,y)-\D_y^{B+1} 
S(x_2,y)=\int_{x_2}^{x_1}\D_y^{B}
F(x,y)\,dx.
\]
Of all the terms arising when I differentiate \rf{factor1} 
$B$ times in $y$, the term in which all derivatives fall on  
$\prod(y-Y_j(x))$ will dominate on the support of $T_x$ after 
a possible reduction of $V$. It follows
that on the support of $T_x$ 
\[ 
|\D_y^{B}F|\approx x^{A+\sum_\ga n_\ga\gg_\ga}=x^{A_1},
\]
where $(A_1,B)$ is the common vertex
of the horizontal infinite edge of $\GG(F)$ and its first compact edge $\ga=1$.

By the previous remarks,
\[
|\GF^{(B+1)}(y)|\gtrsim\bigl|x_1^{A_1+1}-x_2^{A_1+1}\bigr|.
\]
In the case $j'=j$ I have $
\bigl|x_1^{A_1+1}-x_2^{A_1+1}\bigr|\approx 2^{-jA_1}
|x_1-x_2|$ on the support of $K(x_1,x_2)$, so Lemma \ref{vandercorput} 
gives
\[
|K(x_1,x_2)|\<2^{jA_1/(
 B+1)}(\gl
|x_1-x_2|)^{-1/(
 B+1)}.
\]

I apply the following variant of the Schur test (see e.g.
\cite{Halmos}, Theorem 5.2).
\begin{Lem}
\label{schur} 
Let $T$ be an integral operator on $L^2(\rr)$ with kernel $K(x,y)$,
\[
Tf(x)=\int_{-\infty}^\infty K(x,y)f(y)\,dy.
\]
Assume that the quantities
$$
M_1=\sup_y\int|K(x,y)|\,dx,\qquad
M_2=\sup_x\int|K(x,y)|\,dy
$$
are finite. Then $T$ is bounded with $\nt\<(M_1 M_2)^{1/2}$.
\end{Lem}

By this lemma and the estimate of $K(x_1,x_2)$ I have just obtained,    
$$
\|T_j T_j^*\|\lesssim 2^{jA_1/(
 B+1)}\int_0^{2^{-j}}(\gl t)^{-1/(
 B+1)}\,dt
\lesssim \gl^{-1/(B+1)} 2^{j(A_1-B)/(B+1)}.
$$
This of course implies the estimate
\begin{equation}
\label{est1}
\|T_j\|\lesssim  \gl^{-1/(2B+2)} 2^{j(A_1-B)/(2B+2)}.
\end{equation}
As usual, by Lemma \ref{size} I also have a size estimate:
\begin{equation}
\label{est2}
\|T_j\|\lesssim 2^{-j(1+\gg'')/2}\<2^{-j(1+\gg_{1})/2}.
\end{equation}
As the reader may check, 
taking the geometric mean of these two bounds which
kills the $j$-factor gives exactly $\|T_j\|\lesssim
\gl^{-\GD_{1}}$,
with $\GD_{\nu}$ defined as
$$
\GD_\nu=\frac{1+\gg_\nu}{2(1+A_\nu)+2(1+B_\nu)\gg_\nu}.
$$ 
This implies (1) since all $\GD_\nu\>\GD$.

In proving (3), I may assume $j'\>j+2$. Then 
$
\bigl|x_1^{A_1+1}-x_2^{A_1+1}\bigr|\approx 2^{-j(A_1+1)}$ 
on the support of $K(x_1,x_2)$, whence by Lemma \ref{vandercorput}
\[
|K(x_1,x_2)|\lesssim \gl^{-1/(B+1)}2^{j(A_1+1)/(B+1)}=:M.
\]
The support of $K(x_1,x_2)$ is contained in the rectangle of size
$\approx 2^{-j}\times2^{-j'}$. Now Lemma \ref{size} 
gives a bound improved by a factor of $M$: 
\[
\|T_jT_{j'}^*\|\lesssim M2^{-(j+j')/2}=\gl^{-1/(B+1)}2^{j(A_1-B)/(B+1)}
2^{-J/2},
\]
where I denoted $J=j'-j$. By multiplying the estimates
\rf{est2}
for $T_j$ and $T_{j'}$, I get another bound:
\[
\|T_jT_{j'}^*\|\lesssim 2^{-j(1+\gg_{1})}
2^{-J(1+\gg_{1})/2}.
\] 
These two bounds have the form of \rf{est1} and \rf{est2} squared, but
with
an additional factor exponentially decreasing in $J$. Therefore
it is clear that this time taking the geometric mean killing the
$j$-factor
will give 
\[
\|T_jT_{j'}^*\|\lesssim \gl^{-2\GD_{1}}
2^{-\eps J} 
\]  
for some $\eps>0$. This implies (3) and concludes the treatment of $T_x$.


\section{Estimates near $\cal Z$}

This next two sections are devoted to proving upper bounds
for $T^\nu$.

Notice that the sum in \rf{tnu} is almost orthogonal, since the $x$- and
$y$-supports of $T_{jk}$ and $T_{j'k'}$ are disjoint for $|j-j'|$
larger than a fixed constant. Therefore it suffices to
estimate each $T_{jk}$ from the RHS of \rf{tnu} individually.

Fix such a $T_{jk}$. For quite a while the proof is going to proceed
exactly like the argument in the part of Section 2.2.3 dealing with
close rectangles. Analogously to \rf{hessian}, on $R_{jk}$
\beq
\label{fonrjk}
|F|\approx 2^{-jA}2^{-j\gg_\nu B}\prod_{\ga<\nu} 2^{-j\gg_\nu n_\ga}
\prod_{\ga>\nu}
2^{-j\gg_\ga n_\ga}\prod_{i=1}^{n_\nu}|y-Y_{\nu i}(x)| 
\eeq
$$
=2^{-j(\gg_\nu B_\nu+A_\nu-\gg_\nu n_\nu)}\prod_{i=1}^{n_\nu}|y-Y_{\nu
i}(x)|
$$
$$
=2^{-j(\gg_\nu B_\nu+A_\nu-\gg_\nu n_\nu')}\prod_{i=1}^{n_\nu'}|y-Y_{\nu
i}(x)|,
$$
where I ordered $Y_{\nu i}$ so that for $n_\nu'<i\le n_\nu$ we have 
$\Re c_{\nu i}=0$ in $Y_{\nu i}=c_{\nu i}x^{\gg_\nu}+\ldots$.

Let me quickly dispose of the case $n_\nu'=0$, in which I can apply
Lemma \ref{oscest} (the condition (4) is easily checked) and
Lemma \ref{size} to get the oscillatory and size estimates
$$
\|T_{jk}\|\lesssim\gl^{-1/2}2^{j(\gg_\nu B_\nu+A_\nu)/2},
$$
$$
\|T_{jk}\|\lesssim 2^{-j(1+\gg_\nu)/2}.
$$
Now by taking the geometric mean killing the $j$-factor, I obtain the 
required estimate $\|T_{jk}\|\lesssim\gl^{-\GD_\nu}\<\gl^{-\GD}$.  

Now assume that $n_\nu'>0$. Denote $r_i(x)=\Re Y_{\nu i}(x)$, and let
$\bar r _i(x)\in\rr[[x^{1/n!}]]$ be the asymptotic expansion of
$r_i(x)$ at zero. By E.~Borel's theorem, I can find real functions
$f_i(x)$ such that $f_i(x^{n!})\in\cin$ and $f_i(x)\sim \bar r_i(x)$
as $x\to 0$. Moreover, there is one case when I may and will take
simply $f_i(x)=Y_{\nu i}(x)$. Namely, by Proposition 3.1, parts (4),(5),  
this is possible if the series 
$\bar Y_{\nu i}(x)$ is real and different from any other $\bar
Y_{\nu i'}(x)$.

Let $W$ be the union of the graphs of $f_i(x)$ inside $R_{jk}$:
\[
W=\bigcup_{i=1}^{n_\nu'}\bigl\{(x,y)\in R_{jk}\bigl|
y=f_i(x)\bigr\}.
\]
It is not difficult to see that on $R_{jk}$
\[
f_i'(x)\approx x^{\gg_{\nu}-1}\approx2^{-j(\gg_\nu-1)}=:L.
\]
This suggests to consider a Whitney-type decomposition of
$R_{jk}\backslash W$ away from $W$ into rectangles of size $2^{-m}\times
L2^{-m}$. The easiest way to do this is to dilate the set
$R_{jk}\backslash W$ along the $y$-axis $L^{-1}$ times, take the
standard
Whitney decomposition into the dyadic squares away from (the dilation
of) $W$, and
contract
everything to the original scale. As a result, I get a covering
\[
R_{jk}\backslash W\subset \bigcup R_l,\qquad R_l\cap R_{jk}\ne\emptyset,
\]
where $R_l$ are rectangles of size $2^{-m_l}\times
L2^{-m_l}$, $m_l\in\zp$, such that the distance from $R_l$ to $W$ in
the
anisotropic norm $|x|+L^{-1}|y|$ is of the order $2^{-m_l}$.

I claim that the rectangles $R_l$ of \emph{fixed} size form an
almost orthogonal family, i.e.~that for each $R_l$ the number of
rectangles $R_{l'}$ with $m_{l'}=m_l$ such that either the $x$- or the
$y$-projections of $R_l$ and $R_{l'}$ intersect is bounded by a fixed
constant independent of $l$. 

\begin{center}
\scalebox{0.5}{
\rotatebox{-90}{
\includegraphics{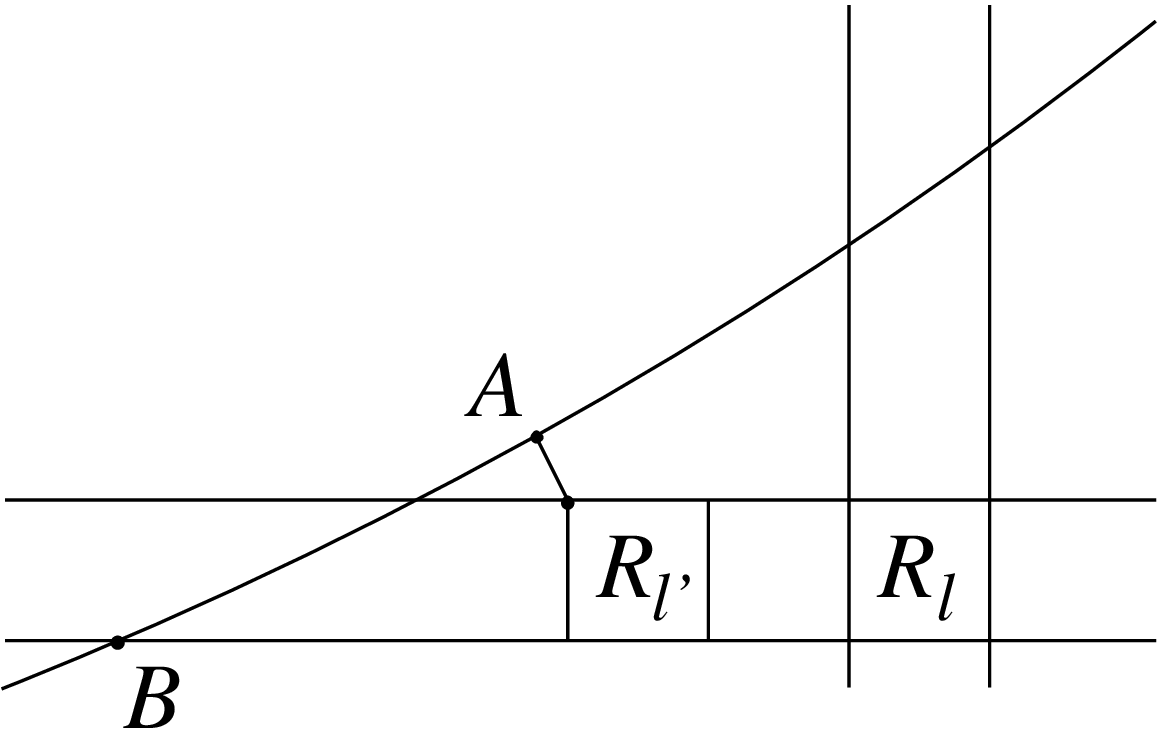}
}}
\\
Fig. 11
\end{center}

Consider the case of intersecting $y$-projections (the other case is
similar). Then $R_{l'}$ is contained in
in the horizontal strip passing through $R_l$ (see Fig.~11).
By dilating along the $y$-axis, I may assume that $L=1$.
Since $\dist(R_{l'},W)\approx 2^{-m_l}$, there exists a point $A$ on
the graph of one of the functions $f_i(x)$ such that
$\dist(R_{l'},A)\approx 2^{-m_l}$. Let $B$ denote the point where the graph
of
$f_i(x)$ intersects the bottom of the strip. Since $f'_i(x)\approx L=1$,
I have $\dist(A,B)\lesssim  2^{-m_l}$, and therefore
$\dist(R_{l'},B)\lesssim 
 2^{-m_l}$. Thus all possible rectangles $R_{l'}$ are situated at a
distance $\lesssim 2^{-m_l}$ from no more than $n_\nu'$ points where the
bottom of the horizontal strip intersects $W$. 
This implies that the number of $R_{l'}$ in the horizontal strip is
bounded by a fixed constant, and the almost orthogonality is verified. 

Now let $R_l^*=(1+\eps)R_l$, where an $\eps>0$ is chosen so small that
$\dist(R_l^*,W)\approx 2^{-m_l}$ (in the anisotropic
norm). Consider a smooth partition of unity $\sum_l\gf_l=1$ on
$\bigcup R_l$ with $\supp \gf_l\subset R_l^*$, satisfying the natural
differential inequalities. I am going to decompose $T_{jk}$ using
this partition of unity. However, this decomposition will not be 
useful near the real multiple branches of $\cal Z$, since I will not have good
control on the size of $F$ there. For now I am just going to
localize away from those branches in the following way.

Let $\gb _i$ denote the power exponent of the first nonzero term
$Cx^\gb$ in the 
asymptotic expansion of $\Im Y_{\nu i}(x)$; $\gb _i:=\infty$ if this
expansion is identically zero. For a large fixed number $Q$ I introduce the set 
$$
W_Q=\mathop{{\bigcup}^*}_{i=1}^{n_\nu'}
\bigl\{(x,y)\in R_{jk}\bigl|
|y-f_i(x)|\<2^{-jQ}\bigr\},
$$
where * indicates that the union is taken over all $i$ such that
$\gb_i=\infty$ and $f_i(\cdot)\ne Y_{\nu i}(\cdot)$. By the choice of
$f_i(x)$, this may happen only if the series $\bar Y_{\nu i}(x)$ is
real and there are several $\bar Y_{\nu i'}(x)$ having $\bar Y_{\nu i}(x)$ as
their asymptotic expansion. One can say that $W_Q$ is a tubular 
neighborhood of width $2^{-jQ}$ of the real multiple branches of $\cal
Z$ (see Fig.~12). 

\begin{center}
\scalebox{0.5}{
\rotatebox{-90}{
\includegraphics{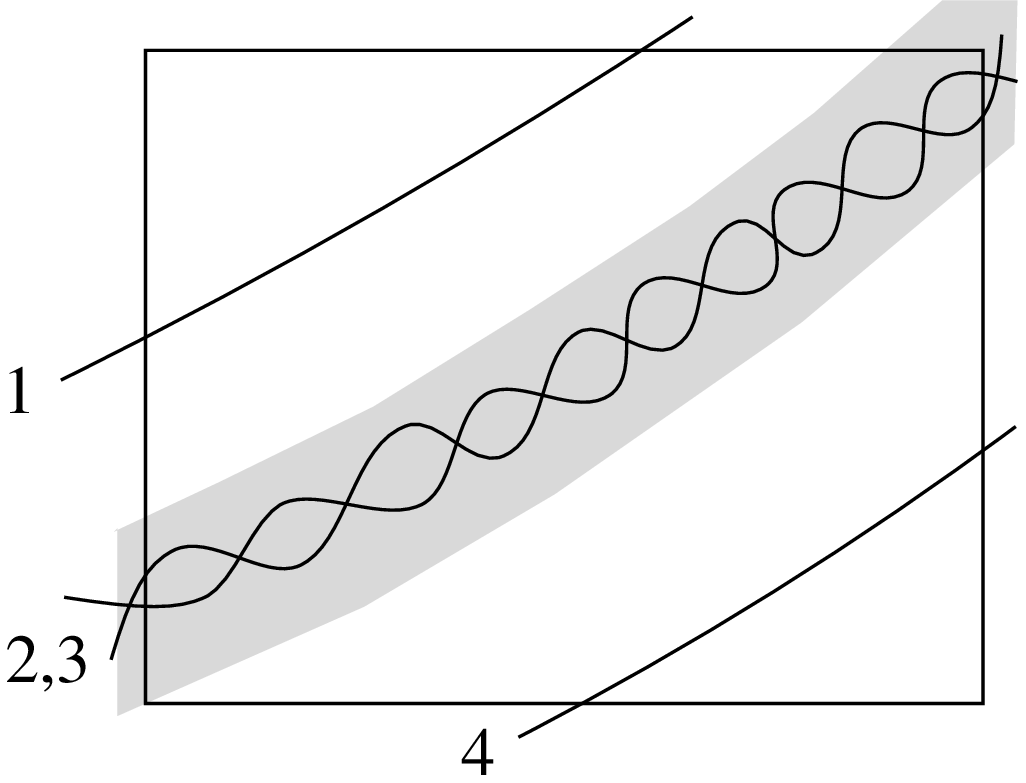}
}}
\\
Fig. 12
\end{center}

The purpose of introducing $W_Q$ is that  
on $R_{jk}\backslash W_Q$ I have (if $j$ is large 
enough, which can be achieved by a further contraction of $V$)
\begin{equation}
\label{yminusfi}
|y-Y_{\nu i}(x)|\approx |y-f_i(x)|+2^{-j\gb_i}. 
\end{equation}  
Now let $\chi_Q$ be a smooth cut-off supported in the double of $W_Q$,
$\chi_Q\equiv 1$ on $W_Q$. I consider the decomposition
\beq
\label{TQ}
T_{jk}=T_Q+T^Q,
\eeq 
$$
T_Qf(x)=\int e^{i\gl
S(x,y)}\chi_Q(x,y)\chi_j(x)\chi_k(y)\chi(x,y)f(y)\,dy,
$$
$$
T^Q=\sum_l T^Q_l,
$$
$$
T^Q_lf(x)=\int e^{i\gl
S(x,y)}\gf_l(x,y)(1-\chi_Q(x,y))\chi_j(x)\chi_k(y)\chi(x,y)f(y)\,dy,
$$
In the rest of this section I prove that
$\|T^Q\|\lesssim\gl^{-\GD}$. The operator $T_Q$ will be dealt with
in the next section.

Let $T_l^Q$ be one of the operators from the decomposition of $T^Q$, and
assume that $T_l^Q\ne 0$, i.e.~that $R_l^*\cap
(R_{jk}\backslash W_Q)\ne\emptyset$.
Fix a point $(x_l,y_l)$ in this last intersection. I claim that 
\begin{equation}
\label{xlyl}
|y-f_i(x)|\approx |y_l-f_i(x_l)|
\end{equation}
for $(x,y)\in R_l^*$ and $i=1,\ldots,n_\nu'$. Indeed, let $(x',y')$
and $(x'',y'')$ be points of $R_l^*$ for which the value of
$|y-f_i(x)|$ is respectively minimal and maximal. Then
$$
|y''-f_i(x'')|\<|y'-f_i(x')|+|y''-y'|+|f_i(x'')-f_i(x')|
$$
$$
\lesssim |y'-f_i(x')|+L2^{-m_l}\lesssim |y'-f_i(x')|,
$$
since $L^{-1}|y'-f_i(x')|\> \dist(R_l^*,W)\gtrsim 2^{-m_l}$. From this
\rf{xlyl} follows.

Now from \rf{fonrjk} and \rf{yminusfi} we see that on $R_l^*$
\[
|F|\approx 2^{-j(\gg_\nu B_\nu+A_\nu-\gg_\nu n_\nu')}
\prod_{i=1}^{n_\nu'}(|y_l-f_i(x_l)|+2^{-j\gb_i}) =:\mu_l.
\]
It follows by Lemma \ref{oscest} (the condition (4) needs to be
checked, but this is easy) that $\|T^Q_l\|\lesssim (\mu_l\gl)^{-1/2}$.

I can get a lower bound on $\mu_l$ by noting that
$|y_l-f_i(x_l)|\gtrsim L2^{-m_l}=2^{-j(\gg_\nu-1)-m_l}$.
This gives
\[
\mu_l\gtrsim 2^{-j(\gg_\nu B_\nu+A_\nu)} 2^{-(m_l-j)n_\nu'},
\]
and therefore
\begin{equation}
\label{tqlosc}
\|T^Q_l\|\lesssim \gl^{-1/2} 2^{j(\gg_\nu B_\nu+A_\nu)/2} 2^{(m_l-j)n_\nu'/2}.
\end{equation}
On the other hand, by Lemma \ref{size},
\begin{equation}
\label{tqlsize}
\|T^Q_l\|\lesssim  2^{-m_l-j(\gg_\nu-1)/2}.
\end{equation}

Now it remains to resum the last two estimates by splitting the family
of operators $T^Q_l$ into almost orthogonal families $l=const$. This
is done exactly how I did it in Section 2.2.3 after Eq. \rf{size22}
\footnote{This resummation was unfortunately done in a wrong way in
my paper \cite{R}. The part of that paper from Eq. (5.6) and until the
end of Section 5 has to be thrown out and substituted by the more
careful argument I give in Section 2.2.3 of this thesis.}.
    
This ends the proof of $\|T^Q\|\lesssim\gl^{-\GD}$.


\section{Estimates near multiple real branches}

To finish the proof of the theorem, I must estimate the operator
$T_Q$ appearing in the decomposition \rf{TQ} of $T_{jk}$. 

In the estimates below I can assume that $\gg_\nu\>1$, since this can
be achieved by passing to the adjoint operator if necessary.

Further, I can assume that $Q$ is chosen so large that the branches
of $\cal Z$ having different asymptotic expansions become completely separated in
the definition of $W_Q$. Since such branches can be treated separately,
I am reduced to the case when $T_Q$ has the form
$$
T_Qf(x)=\int e^{i\gl S(x,y)}\chi_{jkQ}(x,y)f(y)\,dy,
$$
$$
\chi_{jkQ}(x,y)=\chi_j(x)\chi_k(y)\gf(2^{jQ}(y-g(x))).
$$
Here $\gf(t)$ is a $\cin$ cut-off supported in $[-1,1]$,
$g(x^{n!})\in\rr\cin$, $g(x)=cx^{\gg_\nu}+\ldots$, $c\ne0$, and in 
factorization \rf{factor1} exactly $N\>2$ functions $Y_{\nu i}(x)$
have asymptotic expansion coinciding with that of $g(x)$. I will
assume that this happens for $i=1,\ldots,N$. I also re-denote
$W_Q=\{(x,y)\in R_{jk}||y-g(x)|\<2^{-jQ}\}$.

I write $F(x,y)$ as
\[
F(x,y)=\tilde U(x,y)P(x,y),
\]
where $P(x,y)=\prod_{i=1}^N(y-Y_{\nu i}(x))$, and $\tilde U(x,y)$ is
the product of the rest of the terms in \rf{factor1}. 

Since all the branches of $\cal Z$ appearing in $\tilde U(x,y)$ 
are well separated from $W_Q$, there exists a constant $M_1\>0$ such
that
\[
|\tilde U|\approx 2^{-jM_1}\quad\textup{on}\quad W_Q.
\]
Moreover, it can be seen directly that if $\sxy$ is exceptionally
degenerate, we have $M_1=0$.

Further, by Proposition 3.1, parts (4), (5), I know that
$P(x,y)\in\rr\cin_+\(x^{1/n!}\)[y]$, so that $P(x,y)$ is $\cin$ in both
variables on $W_Q$. It is clear that 
\begin{equation}
\label{dyn}
\D_y^N P(x,y)=const\ne0. 
\end{equation}
I claim that, more generally,
\begin{equation}
\label{dxdy}
\D_x^k\D_y^{N-k} P(x,y)\ne 0\quad\textup{on}\quad W_Q,\quad
k=0,\ldots,N.
\end{equation}
Denote $Q(x,y)=P(x^{n!},y)\in\cin_+\(x\)[y]$. The Taylor series of
$Q(x,y)$ is
\[
\bar Q(x,y)=\prod_{i=1}^N(y-\bar G(x)),\qquad \bar G(x)=\bar
g(x^{n!}).
\]
It is clear that
$$
[\D_x^l\D_y^{N-k}\bar Q](x,\bar G(x))=0,\qquad0\<l<k,
$$
$$
[\D_x^k\D_y^{N-k}\bar Q](x,\bar G(x))\ne0.
$$
Therefore the factorizations of $\D_x^l\D_y^{N-k} Q(x,y)$, $l<k$, which can
be obtained as described in the proof of Lemma \ref{pfactor}, will
contain branches with the asymptotic expansion $\bar G(x)$, while the
factorization of $\D_x^k\D_y^{N-k} Q(x,y)$ will not contain such
branches. This implies \rf{dxdy}, provided that $Q$ is large enough,
since $\D_x^k\D_y^{N-k} P(x,y)$ can be expressed as 
\[
(\D_x^k\D_y^{N-k} Q)(x^{1/n!},y)+\sum_{l<k}c_l(x) (\D_x^l\D_y^{N-k}Q)(x^{1/n!},y)
\]
with coefficients $c_l(x)$ growing power-like as $x\to0$.

In addition, the above argument gives an estimate
\begin{equation}
\label{dxn}
\D_x^N P(x,y)\>2^{-jM_2}\quad\textup{on}\quad W_Q,
\end{equation}
for some constant $M_2\>0$; $M_2=0$ if $\sxy$ is exceptionally
degenerate.

Denote $\gs_j(x,y)=\frac 1{j!}\D^j_y P(x,y)$. Consider the
decomposition
$$
T_Q=\sumi_{l=-C}T_l,
$$
$$
T_lf(x)=\int e^{i\gl
S(x,y)}\chi_{jkQ}(x,y)\bar\chi_l(\gs_0(x,y))f(y)\,dy,
$$
where $\bar\chi_l(t)$ is the characteristic function of the set
$2^{-l}\<|t|\<2^{-l+1}$, $C$ is a constant.

I am going to prove the estimates:
\beq
\label{tl1}
\|T_l\|\lesssim 2^{-{l}/{N}+{jM_2}/{2N}},
\eeq
\beq
\label{tl2}
\|T_l\|\lesssim \gl^{-1/2}(\log\gl)^{1/2}2^{l/2}l^{N-1/2}2^{jM_1/2}.
\eeq
The required bound for $T_Q$ can then be derived as follows. 

Consider first the exceptionally degenerate case, when 
$M_1=M_2=0$. I have
\[
\|T_Q\|\lesssim\sumi_{l=0}\min\bigl(2^{-l/N},\gl^{-1/2}2^{l/2}(\log\gl)^{1/2}l^{N-1/2}\bigr).
\]
If it were not for the factor of $(\log\gl)^{1/2} l^{N-1/2}$, the terms in
parentheses would become equal for $l=\frac N{N+2}\log_2 \gl$, and I
would have the best possible estimate $\|T_Q\|\lesssim \gl^{-1/(N+2)}$. In the
present situation I am going to lose something, and to optimize the
loss,
I put $l_0=\frac N{N+2}\log_2 \gl-k\log_2\log_2\gl$ with
indeterminate $k$ and have the estimate\footnote{Here I am being
slightly more careful than in \cite{R} and earn a marginal improvement
in the power of $\log \gl$.}
$$
\|T_Q\|\lesssim 2^{-l_0/N}+\gl^{-1/2}2^{l_0/2}(\log\gl)^{1/2}l_0^{N-1/2}
$$
$$
\lesssim
\gl^{-1/(N+2)}\Bigl[(\log\gl)^{k/N}+(\log\gl)^{N-k/2-1/2}\Bigr].
$$
The optimal value of $k$ is $k=\frac{2N^2-N}{N+2}$, which gives
\[
\|T_Q\|\lesssim \gl^{-\frac1{N+2}}(\log\gl)^{\frac{2N-1}{N+2}}
\]
in complete accordance with what is claimed in the theorem.

Assume now that $\sxy$ is not exceptionally degenerate. In this case the
above argument gives in any case the estimate
\[
\|T_Q\|\<C_\eps 2^{jM}\gl^{-\frac1{N+2}+\eps}
\] 
for any $\eps>0$, with some constant $M$. (I do not pursue the possibility of
obtaining a $\log$ factor here, since as I will see in a moment, what I
have is already good enough.) 

I will need the following more general version of Lemma \ref{size},
which can be obtained immediately from Lemma
\ref{schur}.
\begin{Lem}
\label{Size}
 \tup{Phong and Stein \cite{PSmodels}, Lemma 1.6}
Let $T$ be an integral operator with kernel $K(x,y)$, and assume that\\
\tup{1} $|K(x,y)|\<1$,\\
\tup{2} for each $y$, $K(x,y)$ is supported in an $x$-set of measure
$\<\gd_x$,\\
\tup{3} for each $x$, $K(x,y)$ is supported in a $y$-set of measure
$\<\gd_y$.\\
Then $\nt\<(\gd_x\gd_y)^{1/2}$.
\end{Lem}
By this lemma, I certainly have the estimate
\[
\|T_Q\|\lesssim2^{-jQ/2}.
\]
The idea is that now I can take the geometric mean of the last two
estimates killing the $j$-factor and, if $Q$ is very large, this will
introduce only a very small increase in the exponent of $\gl$,
actually tending to zero as $Q\to\infty$. Thus I have 
\[
\|T_{Q_\eps}\|\<C_\eps\gl^{-\frac1{N+2}+\eps}.
\]  
I am going to show that in the case under consideration
$1/(N+2)>\GD$. 
This allows me to choose and fix $Q$ from the very beginning so large that 
$\|T_{Q}\|\lesssim\gl^{-\GD}$, thus proving the theorem.

I show that in fact $1/(N+2)>\GD_\nu$. Indeed, since I already
have $N$ branches whose expansion starts with $cx^{\gg_\nu}$, I know
that $n_\nu\>N$. Therefore $A_\nu=n_\nu\gg_\nu+A_\nu'\>N\gg_\nu$, and 
\[
2\GD_\nu\<\frac{1+\gg_\nu}{1+N\gg_\nu+\gg_\nu}\<\frac 2{N+2},
\]
since $\gg_\nu\>1$. 
Besides that, the equality holds if and only if $\gg_\nu=1$,
$A_\nu=N$, $B_\nu=0$. But this corresponds exactly to the exceptionally
degenerate case, which is excluded.

I now turn to the proof of the claimed bounds for $T_l$. The proof of
\rf{tl1} is easy and is based on the following well-known 

\begin{Lem}\label{Christ}\tup{Christ \cite{Chr}, Lemma 3.3}
Let $f\in C^N[a,b]$ be such that $f^{(N)}\>\mu>0$ on $[a,b]$. Then for
any $\gg>0$
\[
\bigl|\{x\in [a,b]:|f(x)|\<\gg\}\bigr|\<A_N(\gg/\mu)^{1/N},
\]
where the constant $A_N$ depends only on $N$.
\end{Lem}

By this lemma, in view of \rf{dyn} and \rf{dxn}, the kernel of $T_l$
is supported in a $y$-set of measure $\lesssim 2^{-l/N}$ for each
$x$, and in an $x$-set of measure $\lesssim 2^{-l/N+jM_2/N}$ for each
$y$. Now \rf{tl1} follows by Lemma \ref{size1}.

{\bf Seeger's method.}

The proof of \rf{tl2} constitutes the most intricate part of the whole
argument. It is carried out by
a variation of a method developed in Seeger \cite{See1}, Section 3.
The key idea is to take an additional
dyadic localization in $\gs_j$, $1\<j\<N-1$. Let $l$ be fixed;
all constants below will however be independent of $l$.  Let
$\gg=(\gg_1,\ldots,\gg_{N-1})$ be a vector with integer components
$-C\<\gg_i\<l$, $C$ some constant. Denote
\[
\chi_\gg(x,y)=\chi_{jkQ}(x,y)\bar\chi_l(\gs_0(x,y))\prod_{i=1}^{N-1}
\dbarchi_{\gg_i}(\gs_i(x,y)),
\]
where $\dbarchi_{\gg_i}(t)$ is the characteristic function of the
set $2^{-\gg_i}\<|t|\<2^{-\gg_i+1}$ for $\gg_i<l$, and of the
set $|t|\<2^{-l+1}]$ for $\gg_i=l$.

For an appropriate fixed $C$ I have a decomposition
$$
T_l=\sum_\gg T_\gg,
$$
$$
T_\gg f(x)=\int e^{i\gl S(x,y)}\chi_\gg(x,y)f(y)\,dy.
$$
I am going to prove that for each $\gg$
\begin{equation}
\label{tgamma}
\|T_\gg\|\lesssim \gl^{-1/2}(\log\gl)^{1/2}2^{l/2}l^{1/2}2^{jM_1/2}.
\end{equation}
This will imply \rf{tl2}, since the number of $T_\gg$ in the
decomposition of $T_l$ is $\lesssim l^{N-1}$.
 
The kernel of the operator $T_\gg^*T_\gg$ has the form
\[
K(y_2,y_1)=\int
e^{i\gl[S(x,y_2)-S(x,y_1)]}\chi_\gg(x,y_1)\chi_\gg(x,y_2)\,dx.
\]
Assuming that $y_2>y_1$, and using Taylor's formula in $y$ for
$P(x,y)$, I have 
\beq
\label{[]}
[S(x,y_2)-S(x,y_1)]'_x=\int_{y_1}^{y_2}\tilde
U(x,y)P(x,y)\,dy
\eeq
$$
=\int_{y_1}^{y_2}\tilde U(x,y)\Bigl[\sum_{j=0}^N
\gs_j(x,y_1)(y-y_1)^j\Bigr]\,dy
$$$$
=\sum_{j=0}^N\gs_j(x,y_1)\int_{y_1}^{y_2}\tilde U(x,y)(y-y_1)^j\,dy.
$$
Notice that $\int_{y_1}^{y_2}\tilde U(x,y)(y-y_1)^j\,dy\approx
2^{-jM_1}(y_2-y_1)^{j+1}$. So the RHS of \rf{[]} looks like a
polynomial in $y_2-y_1$ with dyadically restricted coefficients. To
handle such polynomials, I need the following variant of Lemma 3.2
from \cite{See1}. I chose to give a proof, since I have found 
one simpler than in \cite{See1}.

\begin{Lem} 
\label{dyadpol}
For an integer $N\>1$, an integer vector
$r=(r_1,\ldots,r_N)$, $r_i\>0$, and a constant $C>0$ consider the set
${\cal P}={\cal P}(r,C,N)$ of all polynomials of the form
$
P(h)=1+\sum_{i=1}^Na_ih^i
$
with real coefficients $a_i$ satisfying 
$$
|a_i|\in[C^{-1}2^{r_i},C2^{r_i}]\quad if\quad r_i>0,
$$
$$
|a_i|\<C\qquad\qquad\quad\quad{if}\quad r_i=0.
$$
Then there exists a constant $B=B(C,N)$, independent of $r$, and a set
$E\in[0,1]$ of the form 
\begin{equation}
\label{E}
E=[0,2^{\gb_1}]\cup[2^{\ga_2},2^{\gb_2}]\cup\ldots\cup[2^{\ga_s},2^{\gb_s}],
\end{equation}
such that\\ 
\tup{1} $\ga_i$, $\gb_i$ are negative integers, $\gb_1<\ga_2<\gb_2<\ldots<\ga_s<\gb_s\<0$,\\
\tup{2} $s\<B${\rm ;} $\gb_1\>-B\max(r_i)${\rm;}
$\bigl[(1-\gb_s)+\sum_{j=1}^{s-1}(\ga_{j+1}-\gb_j)\bigr]\<B$,\\
\tup{3} $|P(h)|\>B^{-1}$ for $h\in E$ for any $P\in\cal P$.
\end{Lem}
\begin{proof} Put $r_0=0$.
Consider the convex set $\GS$ given as the intersection of the
half-planes lying above the lines $y=r_i+ix$,
$i=0,\ldots,N$. The boundary of $\GS$ consists
of two infinite rays contained in straight lines $y=0$ and $y=r_N+Nx$,
and of some (possibly zero) number of compact segments.

Let $A_i$, $i=1,\ldots,n$ be all the corner points of
the boundary of $\GS$ with the $x$-coordinates
$x_1<x_2<\ldots<x_n$. It is clear that $n\<N$. (In Fig. 12a $N=3$,
$n=3$; in Fig. 12b $N=3$, $n=2$.)

\Picture{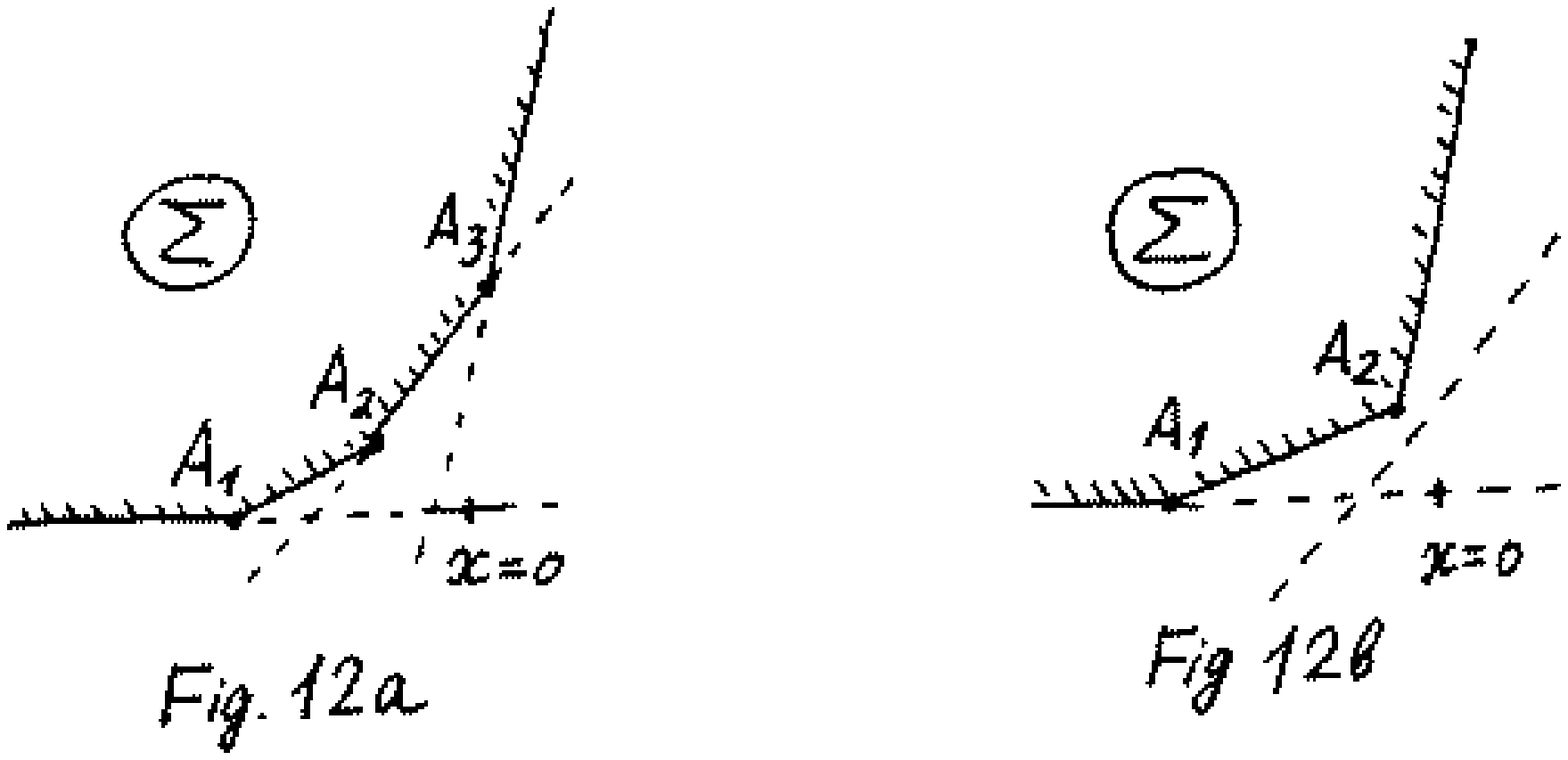}

An observation which will turn out to be important later: if $r_k=0$,
then the line $y=r_k+kx$ cannot contain a compact segment of the
boundary of $\GS$. To see this, it is sufficient to consider how the
lines
$y=kx$, $1\<k\<N-1$, pass with respect to the lines $y=0$ and
$y=r_N+Nx$.
 
I claim that for any $P\in \cal P$ and for large enough $B$
\beql{Star}
|P(h)|\>B^{-1}\quad{if}\quad h\in[0,1],\quad\log_2
 h\notin\bigcup_{j=1}^n(x_j-B,x_j+B).
\eeq

First consider the case 
\beql{lab1}
\log_2 h\in[x_j+B,x_{j+1}-B].
\eeq
Let $k$ be such that the boundary points
$A_j$ and $A_{j+1}$ belong to the line $y=r_k+kx$ (Fig. 12c). 

\Picture{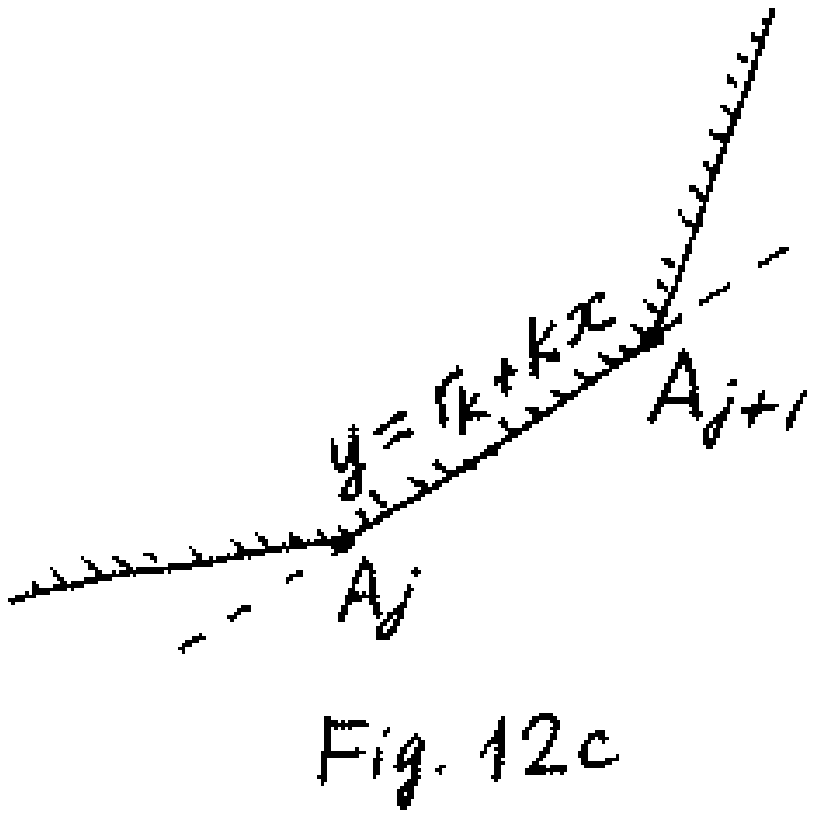}

By the
above observation, $r_k>0$. Since
$A_j$ and $A_{j+1}$ lie above all the other lines $y=r_i+ix$, I have for
all $i$
$$
r_i+ix_j\<r_k+kx_j,\qquad
r_i+ix_{j+1}\<r_k+kx_{j+1}.
$$
From these two estimates it follows that
$$
|a_ih^i|\lesssim |a_kh^k|2^{(i-k)(\log_2 h-x_j)}\qquad(i<k),
$$$$
|a_ih^i|\lesssim |a_kh^k|2^{(k-i)(x_{j+1}-\log_2 h)}\qquad(i>k).
$$
Using \rf{lab1}, I conclude
\[
|a_ih^i|\lesssim |a_kh^k|2^{-|k-i|B}.
\]
This estimate clearly implies
$|P(h)|\gtrsim|a_kh^k|\gtrsim2^{r_k+k\log_2 h}\>1$,
provided that $B$ is large enough.

Second, if $\log_2 h\<x_1-B$, then the same argument as above shows
$|P(h)|\gsim 1$.

Third, if $\log_2 h\>x_n+B$, then if $r_N>0$, I can show in the same
way as above that $|P(h)|\gsim |a_N h^N|\gsim 1$. If $r_N=0$, then
$x_n\>1$, and this region of $h$'s is irrelevant.

So \rf{Star} is verified. Finally, it is not difficult to see that the
exceptional set in \rf{Star} satisfies (1)-(3).
\end{proof}

Now if I take out the factor of
$2^{-l-jM_1}(y_2-y_1)$, the expression in the RHS of \rf{[]} has the form of polynomial in $h=y_2-y_1$
falling under the scope of the lemma with $r_i=l-\gg_i$. So I have a
set $E$ of the form \rf{E} such that
\[
\bigl|[S(x,y_2)-S(x,y_1)]'_x\bigr|\gtrsim 2^{-l-jM_1}(y_2-y_1)\quad
 {if}\quad y_2-y_1\in E.
\]
I claim that this implies
\begin{equation}
\label{|K|}
|K(y_2,y_1)|\lesssim2^{l+jM_1}\gl^{-1}(y_2-y_1)^{-1}\qquad(y_2-y_1\in
 E).
\end{equation}
Indeed, this will follow from Lemma \ref{vandercorput} with $k=1$, if
I prove that there exists a constant $C$ independent of $y_1$ and
$y_2$ such that for fixed $y_1$ and $y_2$\\
(1) the number of intervals of monotonicity of
$[S(x,y_2)-S(x,y_1)]'_x$ considered as a function of $x$ is less than
$C$,\\
(2) the number of intervals comprising the $x$-set where
$\chi_\gg(x,y_1)\chi_\gg(x,y_2)$ is non-zero is less than $C$.

  To show (1), note that $\D_x^N F(x,y)\ne 0$ on $W_Q$. It follows that
\[
\D_x^{N+1}[S(x,y_2)-S(x,y_1)]=\int_{y_1}^{y_2}\D_x^N F(x,y)\,dy\ne 0.
\]
Therefore, $[S(x,y_2)-S(x,y_1)]''_{xx}$ vanishes at most $N-1$ times,
which implies (1).

To show (2), it suffices to check that the number of intervals in the
set $\{x|(x,y)\in W_Q,\ a\<\gs_j(x,y)\<b\}$ is bounded by a constant
independent of $a$ and $b$ for each $0\<j\<N-1$. However, this last
statement follows from \rf{dxdy}.

Unfortunately, to prove the claimed norm estimate for $T_\gg$, I will
need still another decomposition taking into account the form of the
set $E$. Namely, for $1\<k\<s$ and an integer $n$ I put
\[
\chi_{kn}(y)=\gp(2^{\gb_k}y-n), 
\]
where $\gp(t)$ is the characteristic function of the interval $[0,1]$,
and consider the operators 
\[
T_{kn}f(x)=\int e^{i\gl S(x,y)}\chi_\gg(x,y)\chi_{kn}(y)f(y)\,dy.
\]
I am going to prove by induction in $k$ that for each $n$
\[
\|T_{kn}\|\lesssim \gl^{-1/2}(\log\gl)2^{l/2}l^{1/2}2^{jM_1/2}.
\]
The statement for $k=s$ implies the required estimate \rf{tgamma},
since $T_\gg=\sum_n T_{sn}$, and the sum contains no more than
$2^{-\gb_s}\<C$ terms.

For $k=1$, I use the kernel of the operator $T_{1n}^*T_{1n}$, which has the form
\[
\chi_{1n}(y_1)\chi_{1n}(y_2)K(y_2,y_1),
\]
where $K(y_2,y_1)$ is the kernel of $T_\gg^*T_\gg$. If this expression
is not zero, then $|y_2-y_1|\<2^{\gb_1}$. In view of \rf{|K|}, and
also because $|K|\lesssim 1$, Lemma \ref{schur} gives
$$
\|T_{1n}^*T_{1n}\|\lesssim
\int_0^{2^{\gb_1}}\min(1,2^{l+jM_1}\gl^{-1}t^{-1})\,dt
$$
$$
\lesssim2^{l+jM_1}\gl^{-1}\int_0^{\gl2^{\gb_1-l-jM_1}}\min(1,t^{-1})\,dt
\lesssim2^{l+jM_1}\gl^{-1}\log\gl,
$$
which is even better by a factor of $l$ than what I need.

The induction step is performed by using the decomposition
\[
T_{k+1,n}=\sum_{n'}T_{kn'}.
\]
I will need the following variant of the Cotlar--Stein lemma, which
can be proved by an easy adaptation of the standard proof given in
\cite{STEIN}, see e.g.~Comech \cite{Comech}, Appendix.

\begin{Lem} 
\label{cotlarstein}
Let $T_i$ be a family of operators on a Hilbert space
$H$ such that \\
\tup{1} $T_iT_{i'}^*=0$ for $i\ne i'$,\\
\tup{2} $\sum_{i'}\|T_i^*T_{i'}\|\<C$ with a constant $C$ independent of
$i$.\\
Then $\|\sum T_i\|\<C^{1/2}$.
\end{Lem}  

I have $T_{kn'}T_{kn''}^*=0$ for $n'\ne n''$. Let us estimate the sum 
\begin{equation}
\label{sum}
\sum_{n''}\|T_{kn'}^*T_{kn''}\|
\end{equation} 
for a fixed $n'$.  
Since both $T_{kn'}$ and $T_{kn''}$ appear in the decomposition of $T_{k+1,n}$,
I have $|n'-n''|\<2^{\gb_{k+1}-\gb_k}$. Further, the kernel of
$T_{kn'}^*T_{kn''}$ has the form
\[
\chi_{kn'}(y_1)\chi_{kn''}(y_2)K(y_2,y_1).
\]
If this expression is different from zero, then
\begin{equation}
\label{yy}
2^{\gb_k}|y_2-y_1|\in[|n'-n''|-1,|n'-n''|+1].
\end{equation}
Assume first that 
\begin{equation}
\label{nn}
2^{\ga_{k+1}-\gb_{k}}+1\<|n'-n''|\<2^{\gb_{k+1}-\gb_{k}}-1.
\end{equation} 
Then \rf{yy} implies $|y_2-y_1|\in E$, and I can use the estimate
\rf{|K|}. By Lemma \ref{schur},
$$
\|T_{kn'}^*T_{kn''}\|\lesssim\int_{2^{-\gb_k}(|n'-n''|-1)}^{2^{-\gb_k}(|n'-n''|+1)}
2^{l+jM_1}\gl^{-1}t^{-1}\,dt\lesssim 2^{l+jM_1}\gl^{-1}|n'-n''|^{-1}.
$$
Therefore the part of the sum \rf{sum} over $n''$ satisfying \rf{nn}
is bounded by
\[
2^{l+jM_1}\gl^{-1}\sum_{m=2^{\ga_{k+1}-\gb_{k}}}
^{2^{\gb_{k+1}-\gb_{k}}}\frac 1m\lesssim 2^{l+jM_1}\gl^{-1} 
(\gb_{k+1}-\ga_{k+1})\<2^{l+jM_1}\gl^{-1}l,
\]
where I used the fact that by Lemma \ref{dyadpol} (2) $\gb_1\>-B l$.

However, the number of $n''$ which do not satisfy \rf{nn} is bounded
by a constant in view of Lemma \ref{dyadpol} (2), so the
corresponding part of \rf{sum} is bounded by
$C\sup\|T_{kn''}\|^2\lesssim l2^{l+jM_1}(\log\gl)\gl^{-1}$ by the induction 
hypothesis. 

By applying Lemma \ref{cotlarstein}, I complete the induction
step. Theorem 1.1 is now proven.

\chapter{Stopping time}

This chapter stands somewhat separately from the rest of the thesis.
Here I am developing a quite different method of proving upper
norm bounds. This method is incomplete as it stands, and it is unclear
if it is possible to make it complete. In its present form it 
is much less powerful compared to methods based on the geometric analysis 
of the zero set of $\sxy$ which I used above. However, I can use this
method to prove that the estimate \rf{main1} from Theorem 1.1 can be
improved to the optimal $\gl^{-1/4}$ in the case $N=2$.

\section{General idea}

The main idea would be to try to organize an inductive process
which would ``resolve the singularity'' of $\sxy$ by gradually
decreasing the space under its Newton polygon, eventually reducing me
to the non-degenerate case (Fig. 13).

\Picture{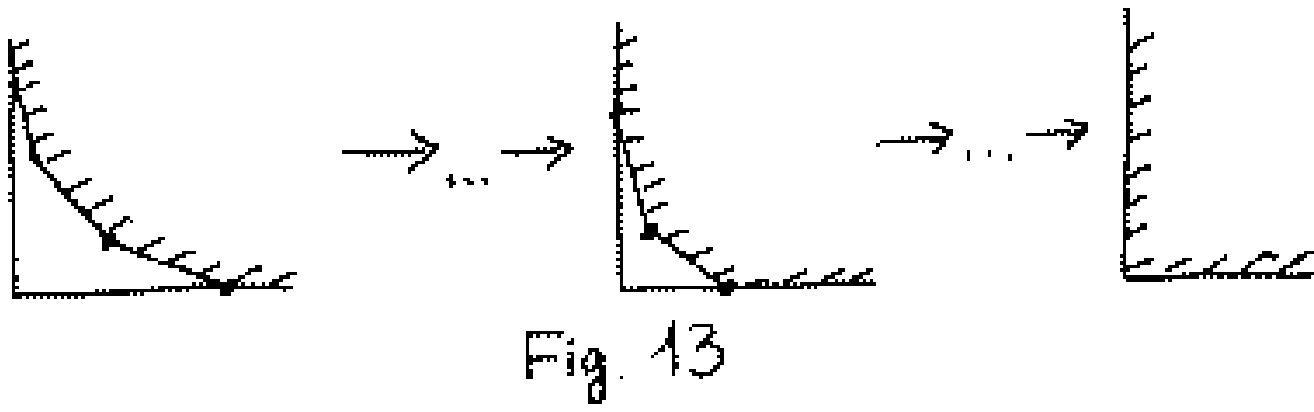}

This idea was first applied to oscilatory integral operators by Phong
and Stein in \cite{PSstop}. Although the proof of that paper is
incomplete as it stands 
(almost orthogonality claim on p.~114 of \cite{PSstop} is
unjustified; see also Remark (c) on p.~150 of \cite{PS}), the argument
can be saved at least in some partial cases \cite{Stein}. Below I 
use a variation of the method 
of \cite{PSstop} and \cite{Stein} to get a somewhat sharper result. 

Still a full realization of the above idea remains elusive. 
The inductive process I
can actually organize works well only for the simplest Newton polygons
consisting of just one edge joining 2 points on the coordinate axes.

Unfortunately, this property may get destroyed already on the first step of
the inductive process (Fig. 14). 

\Picture{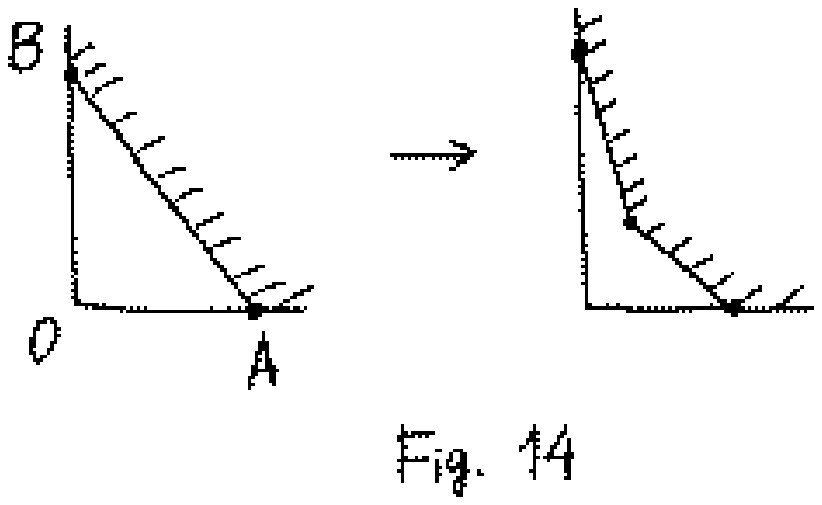}

However, it will not get destroyed provided that {\it there are no integer points lying strictly inside
the triangle $OAB$}. The last condition is satisfied in the following
two cases:
\begi \item $A=1$ or $B=1$ (Fig. 15)
\item $A=B=2$ (Fig. 16)
\ei

\Picture{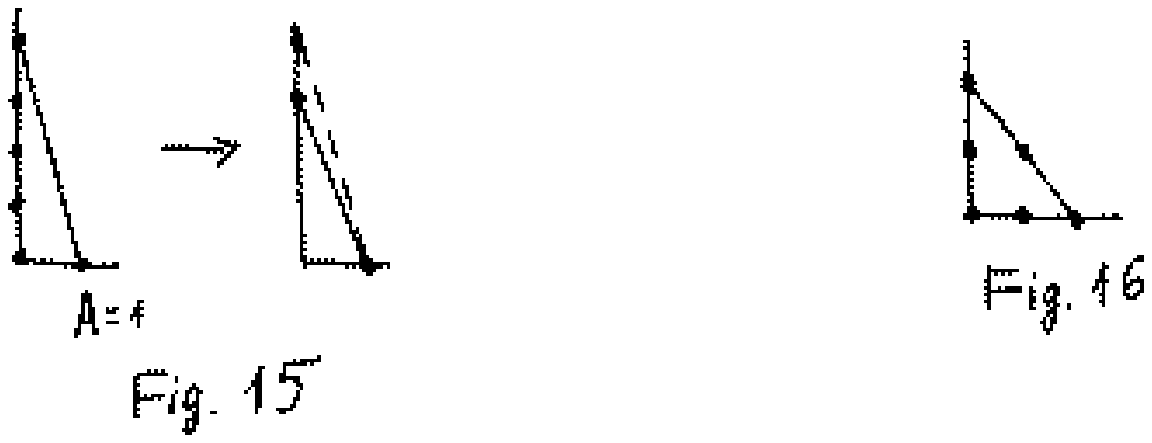}

These are exactly the cases when I am able to produce final results by
this method. In particular, the case $A=B=2$ settles $N=2$ in Theorem
1.1.

\section{Results}

I am going to prove the following

\begin{Th}\label{th1}
  $\ $\\
{\bf I. } Assume that $S(x,y)$ is a smooth phase function in the square
$Q=[0,1]^2$ and that $F=\sxy$ satisfies 
\beql{conditions}
|F^{(1,0)}|\neq 0,\qquad |F^{(0,B)}|\neq 0
\eeq
in Q. Assume also that $\chi$ is a smooth cutoff supported in
$Q$. Then the operator given by \rf{operator} is bounded on $L^2$ with
$$
\ntl\lsim\gl^{-\GD},\qquad \GD=\frac{B+1}{2(2B+1)}.
$$
{\bf II.} If instead of \rf{conditions} I assume that
\beql{20nz}
|F^{(2,0)}|\neq 0,\qquad |F^{(0,2)}|\neq 0
\eeq
in $Q$, then
$$
\ntl\lsim\gl^{-1/4}.
$$ 
\end{Th}

I will need the following somewhat more quantitative auxiliary result,
which implies Part I immediately, and to which Part II will also be
later reduced.

\begin{Th} Let $S(x,y)$ be a smooth phase function in the square
$Q=[0,1]^2$, $\mu$ a real number. Assume that $F=\sxy$ satisfies in $Q$ 
$$
|F^{(1,0)}|\approx \mu,\qquad |F^{(0,B)}|\approx \mu,
$$
and
$$
|F^{(\ga,\gb)}|\lsim\mu
$$
for $(\ga,\gb)\in\{(0,B+1),(0,B+2),(1,1),(1,2),\ldots,(1,B-1)\}.$ 
Assume also that $\chi$ is a smooth cutoff supported in
$Q$. Then the operator given by \rf{operator} is bounded on $L^2$ with
$$
\ntl\lsim(\gl\mu)^{-\GD},\qquad \GD=\frac{B+1}{2(2B+1)}.
$$
\end{Th}

Notice that the number $\GD$ in these results is the Newton decay rate corresponding to the Newton polygon with two
vertices $(1,0)$ and $(0,B)$. Analogously $1/4$ is the right Newton
decay rate for the $(2,0)$--$(0,2)$ Newton polygon. Notice that the $N=2$
exceptionally degenerate phase functions of Theorem 1.1 satisfy
conditions of Theorem 5.2, Part II.

\section{Proofs}

{\bf Proof of Theorem 5.2. } Induction on $B$. For $B=0$ the result
follows from Lemma \ref{oscest}. 

Assume that $B>0$. I divide $Q$ into equal rectangles of size
$\gd^{B}\times \gd$, $\gd=1/2$. If for some of these rectangles the condition $(S)$
below is satisfied, I put it into a numbered collection of rectangles
$\{R_k\}$.
Otherwise I divide it further into equal rectangles of size
$\gd^{B}\times \gd$, now $\gd=1/4$, etc. The stopping condition for a
rectangle $R_k$ of size $\gd_k^{B}\times\gd_k$ is
$$
\min_{\dis R_k^{**}}|F^{(0,n_k)}|\>\mu\gd_k^{B-n_k}\quad\textup{for
SOME}\quad n_k\in\{0,\ldots,B-1\}
\eqno (S)
$$
(star means the doubled rectangle, double star means the quadrupled rectangle).

Eventually all $Q$ up to a set of measure zero becomes
decomposed into rectangles $R_k$. The exceptional set is the
intersection of the zero sets of $F^{(0,n)}$, $n=0,\ldots,B-1$. This set is
of measure zero, since $F^{(0,B)}\neq 0$.

I claim that the covering $\{R_K^*\}$ has finite multiplicity,
that is for every $k$ there are only finitely many $l$'s such that
\beql{nz}
R^*_k\cap R^*_l \neq 0.
\eeq

It is sufficient to prove that \rf{nz} implies $\gd_k\approx\gd_l$. 
Now if $\gd_l\>C\gd_k$, then it follows from 
\rf{nz} that
\beql{diez}
(R^\#_k)^{**}\subset R^*_l,
\eeq
where $R^\#$ denotes the ``parent'' of $R$, that is the rectangle out
of which $R$ was obtained in the $(2^B,2)$-dyadic division process
described above.

But it follows from \rf{diez} that already $R^\#_k$ had to be retained and
not divided further. This contradiction shows that necessarily
$\gd_k\approx\gd_l$,
from which finite multiplicity follows.

Because of finite multiplicity, I can localize the operator
$\tl$ to $R_k^*$ by a smooth partition of unity satisfying the
``right'' differential bounds. Denote the part supported on $R_k^*$ by $T_k$. 

I claim that as well as the lower bound $(S)$, the upper bound
$$
|F^{(0,n)}|\lsim \mu\gd_k^{B-n} \eqno(U_n)
$$
for EACH $n\in\{0,\ldots,B-1\}$ is true on $R_k^*$, and in fact on
$\bar R=(R_k^\#)^{**}$. 

The proof goes like this. Since $R_k^\#$ was not retained, for each
$n=0,\ldots,B-1$ there is a point $(x_n,y_n)$ in $\bar R$
such that
\beql{l*}
|F^{(0,n)}(x_n,y_n)|\<\mu \gd_k^{B-n}.
\eeq
Now by assumption 
\beql{l**}
|F^{(1,n)}|\lsim\mu
\eeq 
in the whole $Q$. Since the $x$-size of $\bar R$ is $\lsim \gd_k^B$,
it follows from \rf{l*} and \rf{l**} by Newton-Leibnitz applied in the
$x$-direction that
\beql{l***}
|F^{(0,n)}(x,y_n)|\<\mu \gd_k^{B-n}\qquad \forall n=0,\ldots,B-1,
\eeq 
provided that $(x,y_n)\in\bar R$. (Fig. 17)

\Picture{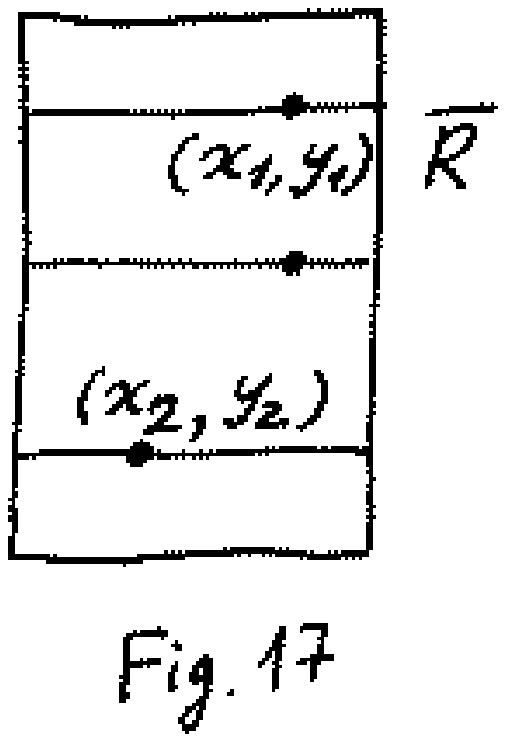}

Now notice that for $n=B$ ($U_n$) is true by assumption in the whole
$Q$,
and that ($U_{n-1}$) follows from ($U_n$) and \rf{l***} by
Newton-Leibnitz
applied in the $y$-direction. So ($U_n$) follows by induction for all
$n$ from $B-1$ to 0.
 
The main reason I need ($U_n$) is to show that for each $\gd$,
the subfamily of rectangles
$R_k^*$ with $\gd_k=\gd$ is
almost orthogonal. 

Indeed, since $|F^{(0,B)}|\>\mu$ on $Q$ and
$|F^{(0,B-1)}|\lsim\mu\gd_k$ on $R_k^*$, by Lemma \ref{Christ} 
there are no more than {\it const} rectangles of the same
$y$-size $\sim \gd_k$ 
with intersecting $x$-projections. 

Analogously, since $|F^{(1,0)}|\>\mu$ on $Q$ and
$|F^{(0,0)}|\lsim\mu\gd_k^B$ on $R_k^*$, by Lemma \ref{Christ} 
there are no more than {\it const} rectangles of the same
$x$-size $\sim \gd_k^B$ 
with intersecting $y$-projections. 

By almost orthogonality, I get
\beql{lsum}
\ntl\lsim\sum_n\sup_{\gd_k=2^{-n}}\|T_k\|.
\eeq
Now the idea is to rescale $T_k$ to a square of size $\sim 1$ by
putting (Fig. 18)
$$
\tilde S(x,y)=S( \gd_k^B x, \gd_k   y).
$$

\Picture{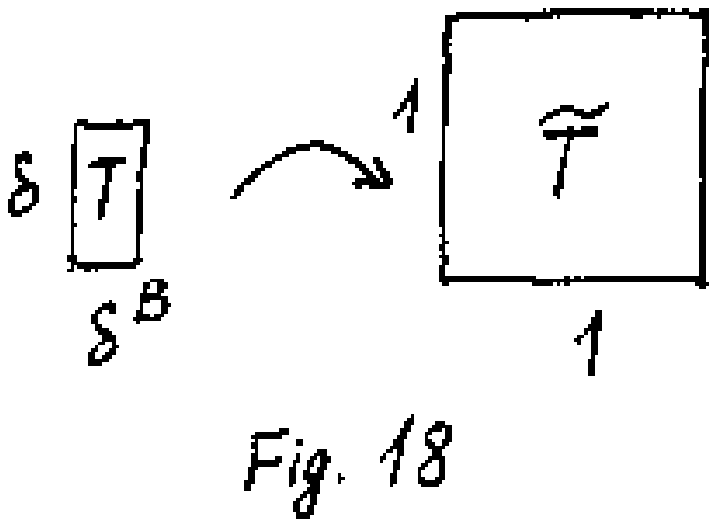}

The norms on $L^2$ are related by
\beql{resc}
\|T_k\|=\gd_k^{(B+1)/2}\|\tilde{T_k}\|.
\eeq
After the rescaling, I get that the phase $\tilde S(x,y)$ satisfies
the conditions of Theorem 5.2 with $n_k$ instead of $B$ and with
$\mu\gd_k^{2B+1}$ instead of $\mu$. Indeed, the main conditions
$$
|\tilde F^{(0,n_k)}|=\gd_k^{B+n_k+1} |F^{(0,n_k)}|\approx\mu\gd_k^{2B+1},
$$ $$
|\tilde F^{(1,0)}|=\gd_k^{2B+1} |F^{(1,0)}|\approx\mu\gd_k^{2B+1}
$$
are satisfied. The auxiliary conditions are checked as follows:

1) $|\tilde F^{(1,\gb)}|=\gd_k^{2B+\gb+1}
|F^{(1,\gb)}|\lsim\mu\gd_k^{2B+1}$ for $\gb=1,\ldots,B-1$, where I used
that by assumption $|F^{(1,\gb)}|\lsim\mu$ on the whole $Q$.

2) $|\tilde F^{(0,n)}|=\gd_k^{B+n+1}
|F^{(0,n)}|\lsim\mu\gd_k^{2B+1}$ for $n=0,\ldots,B$ by ($U_n$), and
for $n=B+1$ by knowing that $|F^{(0,B+1)}|\lsim\mu$ in the whole $Q$.

Thus it follows by the induction hypothesis that
$$
\|\tilde{T_k}\|\lsim\min(1,(\gl\mu\gd_k^{2B+1})^{-\frac {n_k+1}{2(2n_k+1)}}
$$
It follows from \rf{lsum}, \rf{resc}, and the fact that $n_k\<B-1$
that
$$
\ntl\lsim\sum_n 2^{-n(B+1)/2} \min(1,(\gl\mu)^{-\frac B{2(2B-1)}}
2^{\frac {B(2B+1)}{2(2B-1)}n}).
$$
Two progressions balance for 
$$
2^{n_*}\approx (\gl\mu)^{\frac 1 {2B+1}}.
$$
Notice that the second progression is indeed increasing:
$$
\frac {B(2B+1)}{2(2B-1)}-\frac{B+1}{2}=\frac 1{2(2B-1)}>0.
$$
So it follows that
$$
\ntl\lsim 2^{-n_*(B+1)/2}\lsim (\gl\mu)^{-\GD}.
$$
This completes the induction step and the proof of the theorem.
\qedsymbol

{\bf Proof of Theorem 5.1, Part II. } I am going to reduce this result
to the $B=2$ case of Theorem 5.2. This reduction is in fact very
similar to the proof of Theorem 5.2 itself. 

I organize a dyadic decomposition of $Q$, this time into dyadic
squares $R_k$ of size $\gd\times\gd$, $\gd\sim 2^{-n}$, with stopping
condition
$$
\min_{R_k^{**}}|F^{(0,1)}|\>\gd_k\quad\textup{OR}\quad\min_{R_k^{**}}|F^{(1,0)}|\>\gd_k.
\eqno(S')
$$
That is, if $(S')$ is satisfied, the $R_k$ is retained, otherwise it
is further subdivided into 4 squares of equal size, etc.

As before, 
$$
Q=\bigcup R_k
$$
up to a set of measure zero. I show that $R_k^*$ form a covering of
finite multiplicity in the same way as before, and split
$$
\tl=\sum T_k,\qquad \supp T_k\subset R_k^*.
$$

Then I prove that on $R_k^*$
\beql{onr}
|F^{(0,1)}|,|F^{(1,0)}|\lsim\gd_k.
 \eeq
The proof of these bounds is even simpler that that of $(U_n)$. They
 follow immediately by
Newton-Leibnitz from the fact that $R_k^\#$ was not retained.

By Lemma \ref{Christ} I conclude from \rf{onr} and \rf{20nz} that for
each $n$ the $R_k^*$ with $\gd_k=2^{-n}$ form an almost orthogonal
family. This implies \rf{lsum}.

To estimate $\|T_k\|$, I rescale the operator to a square of size
$\sim 1$:
$$
\tilde S(x,y)=S(\gd_k x,\gd_k y),
$$
\beql{102}
\|T_k\|=\gd_k\|\tilde{T_k}\|.
\eeq

Assume that the stopping condition that was actually satisfied for
$R_k$ was $|F^{(1,0)}|\>\gd_k$ (the case of $F^{(0,1)}$ being completely
analogous
because of the $x$-$y$ symmetry). Then after rescaling
$$
|\tilde F^{(1,0)}|=\gd_k^3| F^{(1,0)}|\approx\gd_k^4,
$$$$
|\tilde F^{(0,2)}|=\gd_k^4| F^{(0,2)}|\approx\gd_k^4.
$$

Conditions 
$$
|\tilde F^{(\ga,\gb)}|\lsim\gd_k^4,\qquad
 (\ga,\gb)\in\{(0,3),(0,4),(1,1),(1,2)\},
$$
are also easily checked. So we see that the $\tilde{T_k}$ satisfies
the assumptions of Theorem 5.2 for $B=2$ and $\mu=\gd_k^4$.

It follows that 
$$
\|\tilde{T_k}\|\lsim(\gl\gd_k^4)^{-3/10}.
$$
Going back to \rf{lsum} and \rf{102},
$$
\ntl\<\sum_n 2^{-n}\min(1,\gl^{-3/10}2^{6n/5}).
$$
The progressions are balanced for $2^{n_*}=\gl^{1/4}$, and thus
$$
\ntl\lsim2^{-n_*}=\gl^{-1/4}.
$$
The theorem is proved. \qedsymbol

\bibliographystyle{plain}

\end{document}